\documentclass[a4paper,12pt]{article}
\usepackage{svg}

\usepackage{a4wide}
\usepackage[utf8]{inputenc}
\usepackage{amsmath,amssymb,amsfonts}
\usepackage{theorem}
\usepackage{color,xcolor}

\usepackage{authblk}
\usepackage{xspace}
\usepackage{url}
\usepackage{graphicx}
\usepackage{lipsum} 
\usepackage{enumitem} 
\usepackage{upgreek}
\usepackage{subfigure}
\usepackage{algpseudocode}
\usepackage{algorithm}
\usepackage{longtable}
\usepackage{tabu}
\usepackage{bm}
\usepackage{adjustbox}
\usepackage{bm}
\usepackage{placeins}
\usepackage{float}
\usepackage[normalem]{ulem}
\usepackage[subtle]{savetrees}
\usepackage{booktabs} 
\usepackage{paracol}
\usepackage{multicol}
\usepackage{listings}
\usepackage{cleveref}
\usepackage{float}
\usepackage{tabularx}
\usepackage[frozencache,cachedir=.]{minted}

\theorembodyfont{\normalfont}

\providecommand{\keywords}[1]
{
  \small	\vspace{0.2cm}
  \textbf{{Keywords---}} #1
}

\title{Performance evaluation of mixed-precision Runge--Kutta methods for the solution of partial differential equations}
\author{Ivo Dravins$^1\footnote{ivo.dravins@rub.de}$, Marcel Koch$^2\footnote{marcel.koch@kit.edu}$, Victoria Griehl$^1$, Katharina Kormann$^1$}
\date{ \small{
    $^1$Faculty of Mathematics, Ruhr University Bochum, Bochum, Germany\\%
    $^2$Scientific Computing Center (SCC), Karlsruhe Institute of Technology, Karlsruhe, Germany} 	\\[2ex]%
    {\normalsize \today} 
}

\begin{document}

\setminted{
    fontsize=\small  
}
\newmintinline{cpp}{}
\newmintinline{cuda}{}

\maketitle

\section*{Abstract}
This work focuses on the numerical study of a recently published class of Runge--Kutta methods designed for mixed-precision arithmetic. We employ the methods in solving partial differential equations on modern hardware. In particular we investigate what speedups are achievable by the use of mixed precision and the dependence of the methods \textit{algorithmic compatibility} with the computational hardware. We use state-of-the-art software, utilizing the Ginkgo library, which is designed to incorporate mixed precision arithmetic, and perform numerical tests of 3D problems on both GPU and CPU architectures. We show that significant speedups can be achieved but that performance depends on solver parameters and performance of software kernels. 

\keywords{mixed-precision solvers, mixed-precision time stepping, Runge--Kutta methods, iterative solution methods, preconditioning, GPU computing, stability and truncation}

\section{Introduction}
Mixed-precision algorithms aim at an improved efficiency by the use of low-precision arithmetic where possible---with high-precision arithmetic where necessary to preserve accuracy. While mixed-precision linear algebra has attracted much attention during recent years (see the survey articles \cite{Abdelfattah2021,Higham_Mary_2022} and references therein), the use of mixed-precision for the solution of time-dependent differential equations is still largely unexplored. In the numerical integration of ordinary differential equation literature, some work has been devoted to mixed precision in Runge--Kutta schemes. 
Grant \cite{Grant2022} derived perturbed Runge--Kutta methods that mix stages of different precision, often times with the goal of combining low-precision implicit steps with high-precision explicit steps. For small systems of ordinary differential equations, Burnett, Gottlieb, Grant \& Heryudono \cite{PerformanceEvaluationRK} provided a performance study and Burnett, Gottlieb \& Grant \cite{Burnett2024} a numerical study of stability. Croci \& Rosilho de Souza \cite{Croci2022}, on the other hand, consider explicit stabilized methods based on evaluations of the right-hand side in varying precision for stiff problems. Another example are the mixed-precision exponential integrators proposed by Balos, Roberts \& Gardner \cite{Balos2023}. 

In this article, we consider the application of the mixed-precision Runge--Kutta schemes proposed by Grant \cite{Grant2022} for the time propagation of partial differential equations. In order to consider a performance portable implementation, we built our implementation on the Ginkgo library~\cite{ginkgo}. Ginkgo is a performance-portable high-performance library for sparse linear algebra routines and supports data types in half, single, and double precision as well as linear algebra routines in various precision. Certain operations in Ginkgo are specifically designed to exploit mixed-precision, namely the compressed basis GMRES (CB-GMRES)~\cite{cb-gmres} and the adaptive precision block-Jacobi preconditioner~\cite{anzt2019adaptive}. Additionally, the Ginkgo library offers an accessor class that allows computing arithmetic operations in a given precision without requiring that the operands are stored in that precision. Through this approach, unnecessary copies and conversions can be prevented, as described in the papers~\cite{anzt2020technical,grutzmacher2023using}.
The implementations used for our numerical experiments as well as the data can be retrieved from \url{https://doi.org/10.5281/zenodo.14329619}.

The layout of this contribution is as follows: \Cref{sec:mixedRK} gives an overview of Runge--Kutta methods, followed by aspects of mixed precision computation presenting the methods of interests and outlines some of the challenges faced in this framework. \Cref{sec:solver} describes the solution procedure for the linear systems that arise, including the preconditioner design. \Cref{sec:Ginkgo} gives implementation-specific details related to the software used, with a focus on performance critical computational components of the code. \Cref{sec:results} presents the results of our performance analysis. Conclusions and directions for future work are presented in \Cref{sec:conclusions}.

\section{Runge--Kutta methods and mixed-precision}\label{sec:mixedRK}
Runge--Kutta (RK) methods comprise a widely used class of time-discretization methods. Given an initial value problem of the form
\begin{equation}\label{ode}
    \frac{d u}{d t} = f(t,u), \quad u(t_0) = u_0,
\end{equation}
an RK method approximates the solution at the next time-step by a linear combination of intermediate variables $k_i$, which we refer to as \textit{stages}. Specifically
\begin{equation}\label{RKstep}
   u_{n+1} = u_n + \tau \sum_{i=1}^q b_i k_i, 
\end{equation}
where the stages satisfy
\[
k_i = f\left( t_n + c_i \tau , u_n + \tau \sum_{j=1}^q a_{i,j} k_j\right).  
\]

Any RK method can be fully characterized by two vectors and a matrix, often written in terms of a \textit{Butcher tableau} as illustrated in \Cref{Butchertab}.
\begin{table}[H]
\centering
\caption{The Butcher tableau, $q$ denotes the number of stages.}  \vspace{0.2cm} \label{Butchertab}
\begin{tabular}{l|llll}
$c_1$    & $a_{11}$ &  $a_{12}$    &  $\hdots$        &     $a_{1q}$     \\
$c_2$    & $a_{21}$ & $a_{22}$ &   $\hdots$       &       $a_{2q}$    \\
$\vdots$ & $\vdots$ & $\vdots$ & $\ddots$ &     $\vdots$      \\
$c_q$    & $a_{q1}$ & $a_{q2}$ & $\hdots$ & $a_{qq}$ \\ \hline
         & $b_1$    & $b_2$    & $\hdots$ & $b_q$
\end{tabular}
= \begin{tabular}{l|l}
${\bf c}$&$A$\\\hline & \rule{0pt}{2.5ex} ${\bf b^{{T}}}$
\end{tabular}.
\end{table}
RK methods can be broadly classified by the structure of the RK-matrix $A$: If $A$ has a strictly lower triangular structure the method is explicit, if it is lower-triangular then it is a diagonally-implicit RK (DIRK). In the most general case, when $A$ is full, the method is called a fully implicit RK method. A notable difference between implicit and explicit methods is that implicit methods necessitate the solution of a (possibly non-linear) system at each time-step while explicit methods need only function evaluations. Implicit methods allow for stronger stability properties and higher order (in relation to the number of stages) but entail a significantly higher numerical expense per time-step. 

In the context of HPC applications there have been attempts to mitigate this additional cost by making use of modern hardware through the use of lower-precision computations. The incorporation of lower precision arithmetic, while offering gains in computation efficiency, presents challenges in the loss of accuracy, which counteracts the gains in high-order convergence intrinsic in implicit methods. Several approaches are possible in mitigating this loss of accuracy, e.g., confining the lower precision calculations to interior parts of the solver, i.e., applying a preconditioner in lower precision within a high-precision (outer) Krylov solver. Another approach, which is the focus of this work, is the use of corrector-type methods, where implicit steps are solved in low-precision followed by explicit corrector steps in high precision with the aim to recover accuracy. In this work we study implementations of the latter type of schemes utilizing state-of-the-art software running on both CPU and GPU hardware. Specifically, we focus on a class of methods introduced by Grant in \cite{Grant2022}. 

\subsection{The implicit midpoint method and explicit corrector steps}

To illustrate the incorporation of mixed precision, we consider the implicit midpoint method with correction following \cite{Grant2022}. This entails adding explicit corrector steps to the classical implicit midpoint method. The addition of explicit corrector steps is equivalent to increasing the number of stages of the methods, allowing us to represent the corrector type methods in terms of their Butcher tableaus. In order to differentiate that some steps are evaluated in low precision arithmetic we introduce the superscript $\epsilon$ to denote these low-precision steps. Examples of the implicit midpoint rule with, and without correction steps is given in \Cref{fig_Aq_rat}.

\begin{figure}[H]
\centering
\[
\renewcommand\arraystretch{1.5}
{ \everymath={\displaystyle}
\begin{array}[b]{c|c}
\frac{1}{2}    &  \frac{1}{2}   \\[8pt] \hline
 \rule{0pt}{10pt}  & 1  
\end{array}} \qquad
{ \everymath={\displaystyle}
\begin{array}[b]{l|ll}
\frac{1}{2}^{(\epsilon)} & \frac{1}{2}^{(\epsilon)} &  0 \\[8pt]
\frac{1}{2}     & \frac{1}{2} & 0 \\[8pt] \hline
 \rule{0pt}{10pt}  & 0 & 1
\end{array}} \qquad
{ \everymath={\displaystyle}
\begin{array}[b]{l|lll}
\frac{1}{2}^{(\epsilon)} & \frac{1}{2}^{(\epsilon)} &  0 & \ \ 0 \\[8pt]
\frac{1}{2}     & \frac{1}{2} & 0 & \  \ 0 \\[8pt] 
\frac{1}{2}     & 0 &  \frac{1}{2} &  \  \ 0 \\[8pt] \hline
 \rule{0pt}{10pt}  & 0 & 0 & \ \ 1
\end{array}}
\]
\caption{Examples of Butcher tableau, implicit midpoint rule (left), mixed precision midpoint rule with one correction step (middle) and two correction steps (right). The $(\epsilon)$ superscript denotes steps to be made in low precision.}\label{fig_Aq_rat}
\end{figure}

We recognize that the addition of corrector steps results in a distinct RK method and thus we cannot expect method properties from the base method to automatically hold also for corrected variations, e.g., the A-stability which is known for the standard implicit midpoint rule must be shown to hold also for the corrected case. In the case of implicit midpoint it can be proven that A-stability is preserved for any number of explicit corrector steps---a full proof is given in Appendix A. Writing out explicitly, the method consists of three steps. \textbf{Step 1} is the low-precision solution of the implicit problem:
\begin{align}\label{implicitmid}
    y^{(1)}_{[0]} = u^n + \frac{\tau}{2} f^{\epsilon}(y^{(1)}_{[0]}) .
\end{align}
This is followed by \textbf{step 2}; a sequence of correction steps, essentially a fixed-point iteration, starting from the solution obtained in step 1. 
\begin{align*}
    y^{(1)}_{[1]} &= u^n  + \frac{\tau}{2} f^{}(y^{(1)}_{[0]}) \\
    y^{(1)}_{[2]} &= u^n  + \frac{\tau}{2} f^{}(y^{(1)}_{[1]}) \\
  & \vdots \\
    y^{(1)}_{[p-1]} &= u^n  + \frac{\tau}{2} f^{}(y^{(1)}_{[p-2]}) .
\end{align*}
Finally, \textbf{step 3} is to evaluate the approximated solution at the next time-step using the final stage. 
\begin{align*}
    u^{n+1} = u^n + \tau  f(y^{(1)}_{[p-1]}) .
\end{align*}

The standard implicit midpoint rule is a second order method, from the order conditions (see e.g., \cite{RKorderCond}) we see that the addition of corrector steps preserves second order but does not increase the order. It should be noted, however, as \textbf{step 2} is a fixed point iteration, the addition of corrector steps could be detrimental to convergence. {Consider as an example the case where \eqref{ode} represents a linear PDE, e.g., $f(t,u) = - \Delta u$. Discretizing the spatial operator then leads to the semi-discrete form;} 
\[
\frac{d\bm{u}}{dt} = K \bm{u} \ ; \ K \in \mathbb{R}^{n \times n} .
\]
\textbf{Step 2} will diverge in the limit $p \rightarrow \infty$ when the spectral radius of the iteration matrix exceeds one, e.g., $\rho(\frac{\tau}{2}K)>1$. Hence the addition of more corrector steps cannot be assumed to always result in a better solution, although, method order is retained for any finite number of corrector steps. 

The implicit midpoint rule is useful for illustrating the idea of high precision correction and also visualizes some of the considerations one must take into account when introducing explicit steps to an implicit method---even before considering the effects of low precision arithmetic. In order to achieve methods of higher orders of convergence we move to methods specifically designed for this purpose. 

\subsection{Perturbed Runge--Kutta methods for mixed precision}\label{sec:mixedp_RK}
In \cite{Grant2022}, three different novel RK methods for mixed precision are presented. All three have four stages and are third order methods. We introduce each and outline the method properties. As detailed in \cite{Grant2022}, in describing their properties we need to take into account how the low-precision approximation $f^{\varepsilon}$ relates to $f$. Specifically, if we assume the low-order approximation to be obtained by simply converting e.g., a 64 bit representation to a 32 bit---that is "chopping" the approximation to the desired precision, then the perturbation 
\[
\mathcal{T} = \frac{f-f^{\varepsilon}}{\varepsilon} 
\]
is not Lipschitz continuous, e.g., its derivatives may not exist. This leads to two separate considerations in the characterizing method, when $\mathcal{T}$ is a well-behaved function and when it is not. Here $\varepsilon$ represents a truncation of precision, e.g., we assume that in each implicit step the solution can only be satisfied up to a given tolerance resulting in an error term of order $\mathcal{O}(\varepsilon)$, which is subsequently mitigated by a correction step. {This is the setting in which the methods are derived in \cite{Grant2022}, in our numerical experiments, in addition to employing mixed precision computation, we also rely on iterative solution methods for the implicit steps---this topic is discussed in detail in \Cref{sec:results}. 

\subsubsection{Method 4s3pA} \label{subsec:4s3pA}
The first method, named 4s3pA, is designed to work with badly behaved $\mathcal{T}$, see \cite{Grant2022}, it consists of four stages, two of which entail implicit steps. It is not A-stable---we discuss the method stability properties in \Cref{sec:stability}.  In the case of a well-behaved $\mathcal{T}$ the method error is of the form
\[
E_A = \mathcal{O}(\tau^3) + \mathcal{O}(\epsilon \tau^3) .
\]
We give the combined precision Butcher tableau and method coefficients in \Cref{ButchertabA}. The superscript $\epsilon$ denotes the correction term introduced to mitigate the error introduced by the low precision perturbation as in \cite{Grant2022}. For this and subsequent methods the $c$-values are obtained by row-sums of the coefficient matrix. 
\begin{table}[H]
\centering
\caption{The Butcher tableau and coefficients of the 4s3pA method}  \vspace{0.2cm} \label{ButchertabA}
\begin{tabular}{c|cccc}
$c_1$   & $a_{11}^{\epsilon}$ &  0  & 0  & 0  \\
$c_2$    & $a_{21}$ & 0 & 0 & 0 \\
$c_3$ & $a_{31}+a_{31}^{\epsilon}$ & $a_{32}$ & $a_{33}^{\epsilon}$ &   0     \\
$c_4$    & $a_{41}$ & $a_{42}$ & $a_{43}$ & 0 \\ \hline
         &  0  & 1/2    &  0  & 1/2
\end{tabular} \ \begin{tabular}{c|c}  
$a_{1,1}^{\epsilon} = \ 0.788675134594813 $  & $a_{2,1} =0.211324865405187$  \\ 
$a_{3,1}^{\epsilon} = \ 0.051944240459852 $  & $a_{3,1} =0.709495523817170 $  \\ 
$a_{3,2} = -0.86531425061942 $  & $a_{3,3}^{\epsilon} =0.788675134594813$  \\ 
$a_{4,1} = \ 0.705123240545107 $  & $a_{4,2} = 0.943370088535775$  \\ 
$a_{4,3} = -0.859818194486069 $  &  \\ 
\end{tabular} 
\end{table}
Here we have overlapping high-low precision components, the step-wise implementation is (cf., \cite{PerformanceEvaluationRK}) 
\begin{align}
y^{(1)} &= u^n + \tau a_{1,1}^{\epsilon} f^{\epsilon}(y^{(1)}) \label{mpAstep1} \\
y^{(2)} &= u^n + \tau a_{2,1} f(y^{(1)}) \label{mpAstep2} \\
y^{(3)} &= u^n + \tau [ a_{3,1} f(y^{(1)}) +   a_{3,2} f(y^{(2)})] + \tau [ a_{3,1}^{(\epsilon)} f^{\epsilon}(y^{(1)}) +   a_{3,3}^{(\epsilon)} f^{\epsilon}(y^{(3)})]  \label{mpAstep3} \\
y^{(4)} &= u^n + \tau [ a_{4,1}f(y^{(1)}) + a_{4,2}f(y^{(2)}) + a_{4,3}f(y^{(3)})  ] \label{mpAstep4} \\ 
u^{n+1} &= u^n  + \frac{\tau}{2}(  y^{(2)} + y^{(4)}  )  \label{mpAstep5}
\end{align}

\subsubsection{Method 4s3pB}

The second method, named 4s3pB, is a four-stage, third-order method and the error behaves as
\[
E_B =  \mathcal{O}(\tau^3) + \mathcal{O}(\epsilon \tau^2) ,
\]
that is, we expect a higher contribution from the perturbation error as compared to method 4s3pA. The method is constructed by a perturbation of an L-stable method with the resulting method being A-stable as detailed in \cite{Grant2022}. The Butcher tableau is given in \Cref{ButchertabB}
\begin{table}[H]
\centering
\caption{The Butcher tableau of the 4s3pB method}  \vspace{0.2cm} \label{ButchertabB}
\begin{tabular}{c|cccc}
$c_1$   & $a_{11}^{\epsilon}$ &   0   & 0  &  0     \\
$c_2$   & $a_{21}+a_{21}^{\epsilon}$ & $a_{22}^{\epsilon}$ &  0   &    0    \\
$c_3$   &$a_{31}+a_{31}^{\epsilon}$ & $a_{32}+a_{32}^{\epsilon}$ &  $a_{33}^{\epsilon}$ &  0    \\
$c_4$   & $a_{41}+a_{41}^{\epsilon}$ & $a_{42}+a_{42}^{\epsilon}$ & $a_{43}+a_{43}^{\epsilon}$ &  $a_{44}^{\epsilon}$ \\ \hline
 &    $3/2$  & $-3/2$    &  $1/2$  & $1/2$
\end{tabular} \end{table} \noindent and the method coefficients are given by
\begin{table}[H]
\centering  \begin{tabular}{c|c}  
$a_{1,1}^{\epsilon} = a_{2,2}^{\epsilon} = a_{3,3}^{\epsilon} = a_{4,4}^{\epsilon}=1/2$  &  \\ 
$a_{2,1}^{\epsilon} = -2.376349376129689 $  & $a_{2,1} = 2.543016042796356$  \\ 
$a_{3,1}^{\epsilon} =-2.951484396921318 $  & $a_{3,1} = 2.451484396921318 $  \\ 
$a_{3,2}^{\epsilon} = 0.475891038758779 $  & $a_{3,2} = 0.024108961241221 $  \\ 
$a_{4,1}^{\epsilon} = -0.573861819468268 $  & $a_{4,1} = 2.073861819468268 $  \\ 
$a_{4,2}^{\epsilon} = 0.051944240459852 $  & $a_{4,2} =  2.367724727682735 $  \\ 
$a_{4,3}^{\epsilon} = -1.211868223075524 $  &  $a_{4,3} = 1.711868223075524 .$\\ 
\end{tabular} 
\end{table}
The 4s3pB method has a lower triangular coefficient matrix, i.e., a diagonally implicit Runge--Kutta (DRIK) structure where each stage entails an implicit step. We do not write out the step-wise implementation as the lower/higher precision steps are analogous to method 4s3pA.

\subsubsection{Method 4s3pC}
The final method is the 4s3pC which is derived using simplified order conditions which are less restrictive than the perturbation restrictions used in deriving the A and B methods, see \cite{Grant2022} for details. A resulting property is that in the case of a well-behaved $\mathcal{T}$ the method error is of the form
\[
E_C = \mathcal{O}(\tau^3) + \mathcal{O}(\epsilon \tau^3) ,
\]
while in the case of a badly behaved $\mathcal{T}$, the error is the form 
\[
E_C =  \mathcal{O}(\tau^3) + \mathcal{O}(\epsilon \tau^2) .
\]
The Butcher tableau is given in \Cref{ButchertabC}, as in the case of the 4s3pB method we have a DRIK structure of the coefficient matrix. 
\begin{table}[H]
\centering
\caption{The Butcher tableau of the 4s3pC method}  \vspace{0.2cm} \label{ButchertabC}
\begin{tabular}{c|cccc}
$c_1$    & $a_{11}^{\epsilon}$ &  0    & 0 &   0  \\
$c_2$    & $a_{21}+a_{21}^{\epsilon}$ & $a_{22}^{\epsilon}$ &  0  &  0    \\
$c_3$   &$a_{31}+a_{31}^{\epsilon}$ & $a_{32}+a_{32}^{\epsilon}$ &  $a_{33}^{\epsilon}$ &  0    \\
$c_4$   & $a_{41}+a_{41}^{\epsilon}$ & $a_{42}+a_{42}^{\epsilon}$ & $a_{43}+a_{43}^{\epsilon}$ &   $a_{44}^{\epsilon}$  \\ \hline
         & $b_1 $    & $b_2$    &  $b_3$  & $b_4$
\end{tabular}
\end{table}
\noindent The coefficients are given by 
\begin{table}[H]
\centering \begin{tabular}{l|l|l}  
$a_{1,1}^{\epsilon} = 0.511243008730995 $ &   &   $ b_1 = 0.002837446974069$      \\ 
$a_{2,1}^{\epsilon} = -1.999347282862640 $  & $a_{2,1} = -0.050470366527530$   & $ b_2 = 0.336264433650450 $ \\ 
$a_{2,2}^{\epsilon} = 1.957161067302390  $   &     &   $ b_3 = 0.806376720267787 $ \\
$a_{3,1}^{\epsilon} = 0.443312893511937 $  & $a_{3,1} = 0.368613367355336 $  & $ b_4 = -0.145478600892306 $ \\ 
$a_{3,2}^{\epsilon} = -0.573131033672219 $  & $a_{3,2} = 0.273504374252976 $ &   \\ 
$a_{3,3}^{\epsilon} = 0.128283796414019 $  &   &  \\ 
$a_{4,1}^{\epsilon} = -2 $  & $a_{4,1} = 1.803794668975043 $ &  \\ 
$a_{4,2}^{\epsilon} = -0.160330320741428 $  & $a_{4,2} =  0.097485042980759 $ &  \\ 
$a_{4,3}^{\epsilon} = 0.579597314161362 $  &  $a_{4,3} = -1.895660952342050  $ & \\ 
$a_{4,4}^{\epsilon} = 1.484688928981990$  &    & 
\end{tabular} 
\end{table}
 
\subsection{Method stability in the context of mixed precision}\label{sec:stability}

The topic of mathematical and numerical stability in the context of mixed precision arithmetic is an intricate topic which this work does not attempt to explore in-depth. We only outline some central issues. The central motivation of the use of mixed precision arithmetic is the potential speedup gained from working with numbers utilizing fewer bits for storage. This leads us to the first notable challenge, namely overflow and underflow when working with low precision. The largest/smallest normal numbers in the IEEE 754 standard are given in \Cref{tab:normalnumbers}.

\begin{table}[H]
    \centering  
    \caption{Smallest and largest normal number (IEEE 754 binary format)} \vspace{0.2cm}
    \begin{tabular}{l|c|c} 
     \hline 
       Type  &  Smallest normal number & Largest normal number \\ \hline 
       Half (16 bit)  &  $2^{-14} \approx 6.1 \cdot 10^{-5}$  &  $ 2^{15} \cdot (2-2^{1-11}) = 65504 $ \\ 
       Float (32 bit)  & $ 2^{-126}   \approx 1.2 \cdot 10^{-38} $ & $2^{127} \cdot (2-2^{1-24}) \approx 3.4 \cdot 10^{38}$  \\ 
       Double (64 bit)  & $ 2^{-1022} \approx 2.2 \cdot 10^{-308} $ &  $2^{1023} \cdot (2-2^{1-53}) \approx 1.8 \cdot 10^{308} $ \\ 
    \end{tabular}
        \label{tab:normalnumbers}
\end{table}

Our framework supports computations in half precision but as can be read from \Cref{tab:normalnumbers}, the largest normal number in 16 bit precision being 65,504 is very restrictive when solving problems of practical interest, e.g., a spatial discretization parameter $h=0.003$ already entails that the common quantity of interest $1/h^2 > 65504$, leading to failure due to overflow. Similarly, the smallest normal number being of the order $\sim 10^{-5}$, which is larger than e.g., commonly used tolerances, leads to underflow issues.
This issue can be mitigated by re-scaling the problem, trying to keep relevant quantities within an allowable interval, as presented in \cite{higham2019squeezing}, but this strategy introduces an unavoidable overhead which counteracts the speedup of the individual computations. 
In \cite{TurekHalf2022} the authors offer an approach to handle half precision for FEM problems by choosing a specific finite element basis, which is not, however, applicable in our setting.
Thus, we have not developed competitive general implementations for half precision, and focus our analysis on the setting double/float. 

The discussion of method stability is also involved as it entails understanding method sensitivity to round-off and propagation of rounding errors throughout the solution procedure. As we focus on numerical aspects, we will not treat this subject in depth, for a study we refer to \cite{StabilityEvaluationRK}. A simple example, while not capturing all complexity, is, however, useful for the purposes of illustration. 

Consider the method 4s3pA presented in \Cref{subsec:4s3pA}, we separate the coefficient matrix according the high- and low-precision parts $A = A^{(h)}+A^{(\epsilon)}$. We construct the stability function 
\begin{equation}
    R(\mu) = 1 + \mu \bm{b}^{\top} ( \mathbb{I} - \mu A)^{-1} \bm{e},
\end{equation}
where $\bm{e}=[1,...,1]^{\top}$. We numerically evaluate the region of stability, i.e., $R_S :=\{z \in \mathbb{C} : |R(z)| \leq 1 \} $. We do this twice, in the first instance we represent both $A^{(h)}$ and $A^{(\epsilon)}$ in double precision, in the second case we first truncate the coefficients in $A^{(\epsilon)}$ to half precision, followed by an evaluation of $R_S$ fully in double precision. Such an approach clearly does not capture the complexities of the propagation of rounding errors, it does, however, illustrate how significantly the stability region can be affected by the truncation of method coefficients alone. The two stability regions are shown in \Cref{fig:StabilityMPA} where we clearly see that in the case where we truncate the coefficients of $A^{(\epsilon)}$ the stability regions shrinks significantly.

\begin{figure}[H]
\centering
\includegraphics[width=0.49\textwidth]{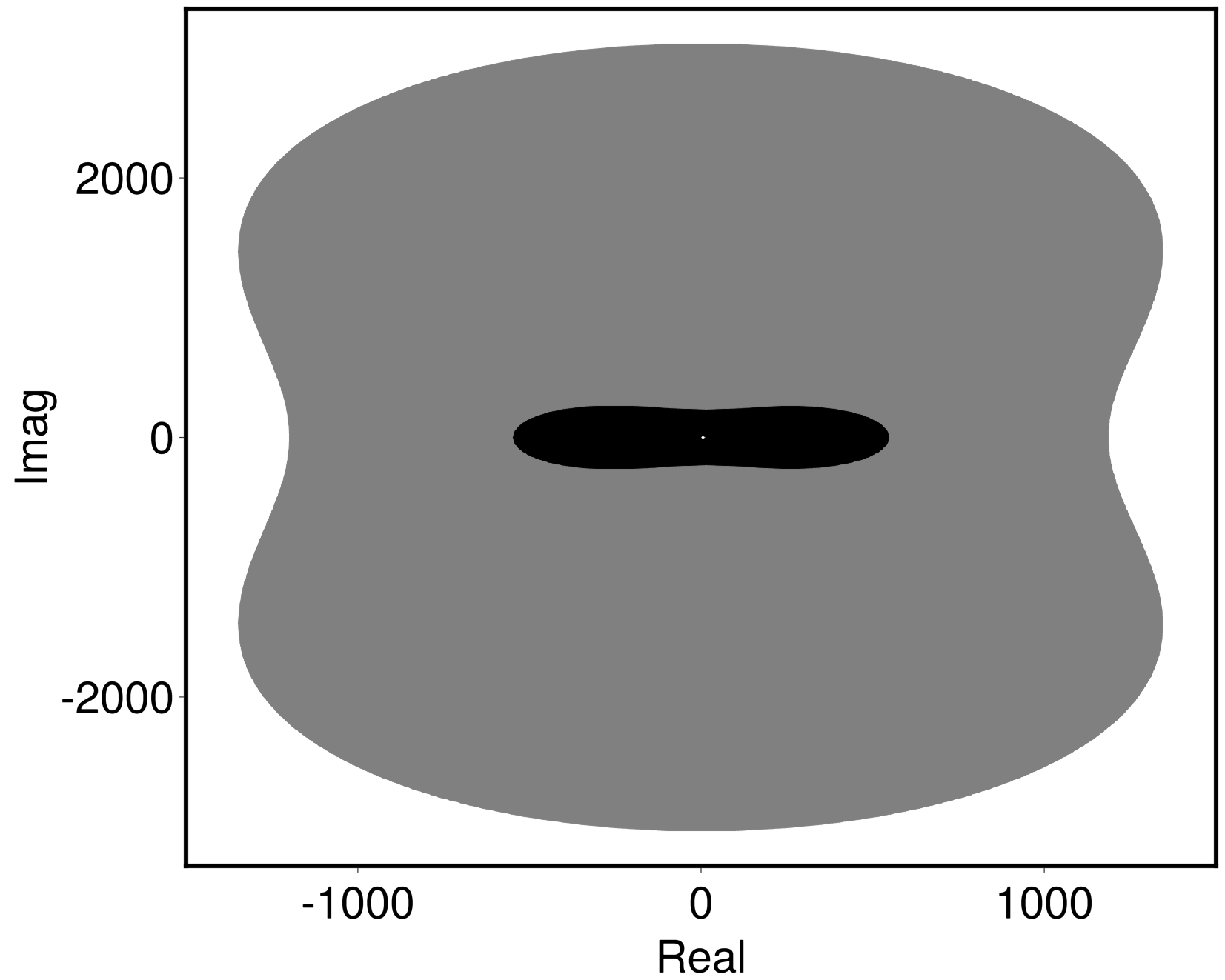}
\caption{Stability regions for the 4s3pA method: 64bit $A^{(\epsilon)}$ - gray, 16bit $A^{(\epsilon)}$ - black.}
\label{fig:StabilityMPA}
\end{figure}

The use of low precision arithmetic for some stages must thus be assumed to fundamentally alter method properties, in addition to introducing larger rounding errors. Note that we illustrate the difference between double and half as the large difference illustrates this clearly, in the case of truncating $A^{(\epsilon)}$ to float precision, the difference in this example is barely perceptible. 

The stability of the B and C methods are insensitive to the truncations of the type described above. Studying the stability region numerically with Julia or Mathematica indicates that both methods are A-stable. The stability regions close to zero are visualized in \Cref{fig:StabilityMPBC}. 
\begin{figure}[H]
\centering
\subfigure[4s3pB]{\includegraphics[width=0.49\textwidth]{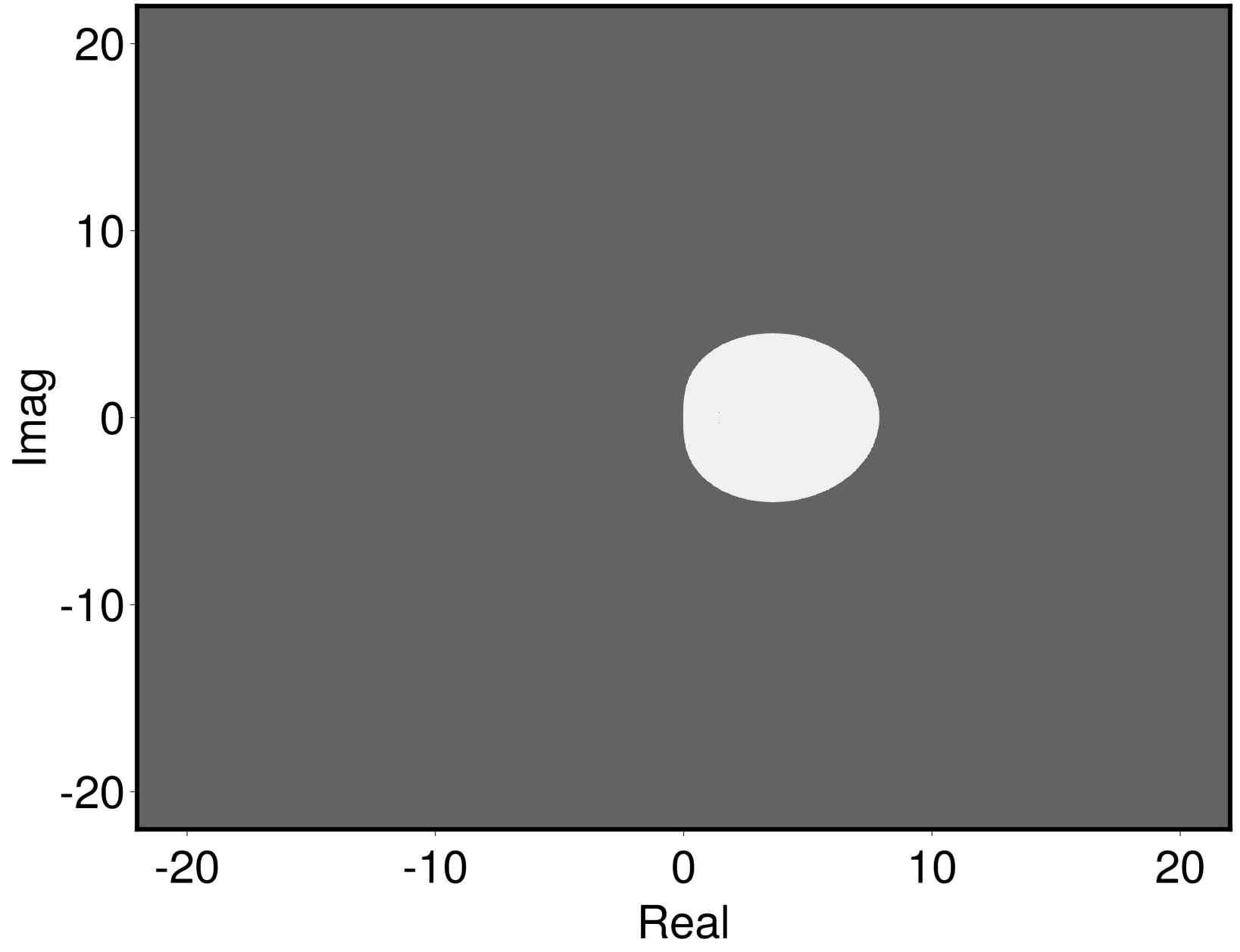}}
\subfigure[4s3pC]{\includegraphics[width=0.49\textwidth]{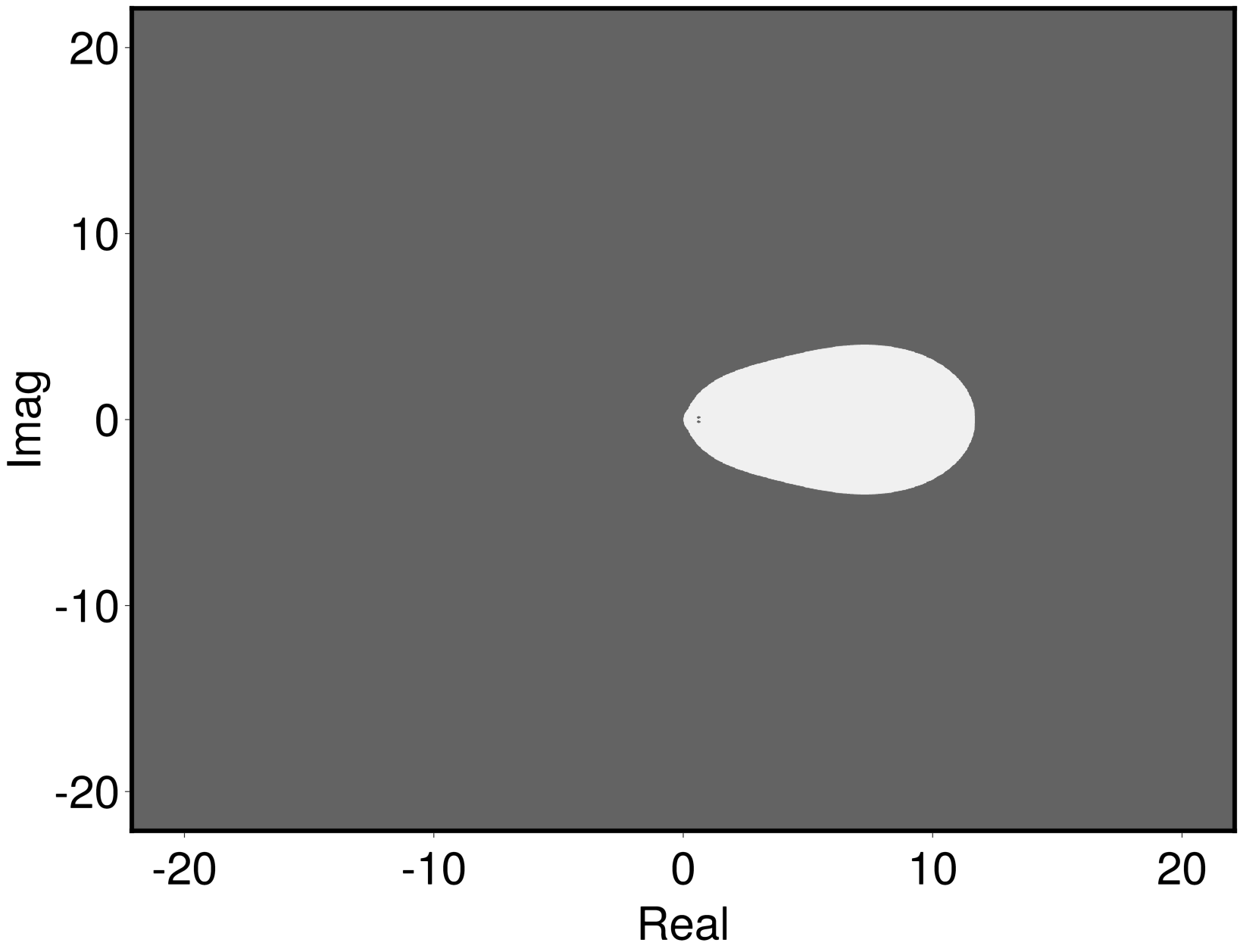}}
\caption{Stability regions for the B and C methods in gray.}
\label{fig:StabilityMPBC}
\end{figure}

\section{Problem setting and solver design} \label{sec:solver}
The core aim of this work is to experimentally investigate performance of the mixed-precision approaches utilizing modern software purpose-built for this application, running on different hardware, e.g., GPUs and CPUs. This isolated aim would not be aided by a non-linear problem which introduces significant algorithmic complications, in particular in the handling of mixed precision in the non-linear solver. In the stability study \cite{StabilityEvaluationRK} non-linear problems are considered and the stopping criteria of the non-linear solver is tied to the precision used, leading to $\sim$15x reductions in runtime---well in excess of 2x speedup expected by use of lower precision alone. We instead opt to focus on two distinct, representative linear problems in 3D. The first being the heat equation, the other being a pure advection example. This setting allows us to explore large scale behavior, giving insight also relevant for the non-linear case. In our single-GPU setting a core constraint of the focus on 3D becomes that spatial refinement is constrained by the GPU memory.    

In the cases considered the semi-discretized form of \eqref{ode} reads 
\[
\frac{d\bm{u}}{dt} = K \bm{u}
\]
where $K \in \mathbb{R}^{N \times N}$. All methods presented in \Cref{sec:mixedp_RK} share similar numerical components, chiefly, in each implicit stage, e.g., \eqref{implicitmid} and \eqref{mpAstep1}, we need to solve a system of the form 
\begin{equation}\label{corelinear}
(a_{ij} \mathbb{I} + \tau K ) \bm{x} = \bm{b} .
\end{equation}
We expect this step to constitute the main numerical cost. As the focus is on large 3D problems we employ an iterative Krylov method in solving \eqref{corelinear}, necessitating also the use of efficient preconditioners. In the context of a performance study for mixed-precision the computationally significant components are critical in enabling speedup, specifically, the operations suitability for particular hardware. In order to quantify this effect, and to have comparable preconditioners for both reference problems, we elect to discretize space with a uniform Cartesian finite difference method (FDM). This allows us to either factorize exactly, or closely approximate the system matrix in \eqref{corelinear} by factorized tensor identities, enabling the use of the so-called fast-diagonalization-method. 

Specifically, as we are solving problems on the unit cube discretized by a uniform Cartesian mesh, we choose a preconditioner of the following form
\begin{equation}\label{preconbasic}
    P = \mathbb{I}_n \otimes \mathbb{I}_n \otimes A + \mathbb{I}_n \otimes B \otimes \mathbb{I}_n  + C \otimes \mathbb{I}_n  \otimes \mathbb{I}_n  ,
\end{equation}
where $A$, $B$ and $C$ are $n \times n$ matrices which represent one dimensional operators in the $x$, $y$, and $z$ directions respectively. This form allows us to use the spectral decomposition (e.g., $A= Q_A \Lambda_A Q_A^{-1}$) to write the inverse action of the preconditioner as
\begin{align}\label{Pinvform}
           P^{-1} &= (Q_C \otimes Q_B \otimes Q_A) \underbrace{(  \mathbb{I}_n \otimes  \mathbb{I}_n \otimes \Lambda_A  +  \mathbb{I}_n \otimes \Lambda_B \otimes  \mathbb{I}_n + \Lambda_C \otimes  \mathbb{I}_n \otimes  \mathbb{I}_n )^{-1}}_{P_D^{-1}} (  Q^{-1}_C \otimes  Q^{-1}_B \otimes  Q^{-1}_A), 
\end{align}
where $P_D$ is a diagonal matrix. We note the following relation
\begin{equation}\label{TensorDecomp}
     ({Q}_C \otimes {Q}_B \otimes {Q}_A) = \underbrace{({Q}_C \otimes I_n \otimes I_n)}_{T_L} \underbrace{(I_n \otimes {Q}_B \otimes I_n )}_{T_M} \underbrace{( I_n \otimes I_n \otimes {Q}_A)}_{T_R}. 
\end{equation}
Thus, in order to calculate the action of $P^{-1}$, we need the action to invert the diagonal matrix $P_D$ and apply operations of the type $T_{L,M,R}$. The three latter tensor operations are especially useful in illustrating factors that affect speedup on different hardware. This is due to the distinct access patterns required by the three different operations as we discuss in depth in the subsequent sections.  

Another benefit in the context of the current study is that the two model problems, while sharing a similar preconditioning approach, result in different computational costs. In the case of the heat equation, the 1-dimensional matrices, $A, \ B$ and $C$ are symmetric positive definite (SPD) resulting in a real spectral decomposition and thus all factors in \eqref{Pinvform} being real. In the convection case we have a complex-spectrum leading to a higher computational cost per application. Furthermore the convection equation necessitates a Krylov solver which allows for non-symmetric matrices, we use the generalized minimal residual method (GMRES), while for the heat equation we employ the computationally cheaper conjugate gradient (CG) to solve the resulting SPD systems. 

\section{Implementation in Ginkgo}\label{sec:Ginkgo}
The following section discusses the implementation of the Runge--Kutta methods using the Ginkgo library~\cite{ginkgo}.
Ginkgo is written in C++ and uses modern software design principles to provide a high-performance, flexible, and sustainable library for sparse linear algebra routines.
It supports all major GPU vendors, i.e., NVIDIA, AMD, and Intel, as well as modern multicore CPUs through OpenMP.

\subsection{Short overview of Ginkgo}

The design and functionality of Ginkgo make it especially useful for the approach discussed in this paper.
Besides efficient hardware-specific kernels, Ginkgo has native support for mixed-precision operations.
Together, these qualities allow investigating the effect of mixed-precision with little overhead.
In addition, the flexible design allows to easily use and combine the components provided by the library, and also provide user-defined operations that can directly interact with the rest of the library.
The implementation of the preconditioner mentioned in \cref{sec:solver} is an example of such an operation.
More details on the handling of mixed-precision and user-defined operations is provided below.

The central concept in Ginkgo to allow user-defined operations is the \emph{linear operator} and the corresponding interface class \cppinline{gko::LinOp}.
It provides an abstraction for mathematical linear operators.
All Ginkgo matrix types are derived from this common type, as well as preconditioners, solvers, reorderings, even though the latter ones might not fit the mathematical definition.
With this abstraction in place, many operations, iterative Krylov solvers included, can be formulated just through the common base class.
Users can provide their own types derived from \cppinline{gko::LinOp} which implement a specialized operator application, in this instance the preconditioner from \cref{sec:solver}.
The custom preconditioner can then be directly used in Ginkgo's Krylov solvers.

The Ginkgo library provides a flexible approach to mixed precision.
First, the precision of most types can be chosen individually, which means a single application can use multiple objects stored in different precisions.
Users can choose \cppinline{gko::half}\footnote{Ginkgo's implementation of the FP16 precision, available since version 1.9.0.}, \cppinline{float}, \cppinline{double}, or their complex counterparts, as precision for their objects.
Additionally, Ginkgo objects of different precision can seamlessly interact with each other.
If the precision of the objects does not match, then the non-matching ones are converted to the precision of the object that an operation is called on.
For example, the vectors in a matrix-vector product will get converted to the precision of the matrix.
Obviously, this requires an additional copy.
Ginkgo provides a mode that circumvents the copy, by using an accessor class to directly access data in different precisions, see~\cite{anzt2020technical,grutzmacher2023using} for more details.
However, it should be noted that this is only available for the matrix-vector product with the ELL or CSR format, and thus could not be leveraged in this work.

\subsection{Implementation of the Mixed-precision Runge--Kutta methods}

In the following, the specifics of implementing the mixed-precision RK methods, introduced in~\cref{sec:mixedRK}, with Ginkgo are discussed.
As mentioned before, Ginkgo allows using different precision for the explicit and implicit parts of the RK methods.
All explicit operations, including linear combination of vectors, are done in high precision, i.e. \cppinline{double} precision, and the stage vectors themselves are also stored in high precision.
Solving the implicit systems can be done in high or low precision.
If it is solved in low precision, i.e. \cppinline{float} precision, then Ginkgo's internal conversion is used to convert the input and output vectors automatically to lower precision and all solver internal operations are executed in this precision.
The output is again converted back to high precision at the end.

The preconditioner mentioned in \cref{sec:solver} is implemented as a custom Ginkgo linear operator.
The tensor-product applications are implemented as independent kernels for each tensor $T_{L,M,R}$, which can be written as
\begin{equation*}
    \left(T_R x\right)_{ijk} = \sum_{q=0}^{n-1} Q_{iq} x_{qjk}, \
    \left(T_M x\right)_{ijk} = \sum_{q=0}^{n-1} Q_{jq} x_{iqk}, \
    \left(T_L x\right)_{ijk} = \sum_{q=0}^{n-1} Q_{kq} x_{ijq}, 
\end{equation*}
where 0-based indexing is used with $i,j,k \in [0, n)$ and with lexicographical numbering of global vectors, i.e., $x_{ijk} = \bm{x}_{i +  jn + kn^2}$.
Thus, the kernels are variations of matrix-vector products, with different strides.
Due to their language similarities, CUDA and HIP kernels are implemented using the same sources, with only minor adjustments for the actually used GPU backend.
The OpenMP kernels are provided separately.
Although Ginkgo also allows for SYCL kernels, those were not implemented as Intel GPUs were not available to either test or benchmark on.

The CUDA/HIP kernel\footnote{Only CUDA terminology is used throught this section.} for $T_M$ are shown in~\Cref{lst:gpu-tpm}. 
The $T_{L, R}$ kernels are nearly identical, only the access to \cppinline{B} and \cppinline{in} differs.
Both kernel uses one thread per row of the in- and output vectors.
In case of $T_R$ this would lead to inefficient memory access patterns, as accessing the input vector would be non-coalesced.
To remedy this issue, a full warp per row is used together with a warp-wise reduction.

\begin{listing}[tbp]
\begin{minted}{cpp}
template<typename ValueType>
__global__ void tensor_product_middle_kernel(
    gko::size_type num_rows, const ValueType* in, ValueType* out,
    const ValueType* B) {
  auto global_index = blockIdx.x * blockDim.x + threadIdx.x;
  auto n = num_rows * num_rows * num_rows;
  auto num_cols = num_rows;
  if (global_index >= n) {
    return;
  }
  // extract tensor indexing
  auto k = global_index / (num_rows * num_rows);
  auto j = (global_index - k * num_rows * num_rows) / num_rows;
  auto i = (global_index - k * num_rows * num_rows) % num_rows;
  auto vector_start = k * num_rows * num_rows + i;
  ValueType acc = 0;
  for (gko::size_type q = 0; q < num_rows; q++) {
    auto vector_index = vector_start + q * num_rows;
    acc = B[j * num_cols + q] * in[vector_index] + acc;
  }
  out[global_index] = acc;
}
\end{minted}
\caption{CUDA/HIP kernel for applying $T_M$.}
\label{lst:gpu-tpm}
\end{listing}

An alternative implementation might use BLAS calls.
Since the tensor product applications can be written as blocked matrix-vector products with strides, the batched interface of the cuBLAS\footnote{\url{https://docs.nvidia.com/cuda/cublas/}}, or rocBLAS\footnote{\url{https://rocm.docs.amd.com/projects/rocBLAS/}} libraries could be considered.
A preliminary exploration of these implementations showed, however, worse performance than the above-mentioned kernels.
Especially the non-unit-stride versions performed poorly.
Thus, they were not considered further for the GPU implementation.

\Cref{lst:cpu-tpm} shows the OpenMP kernel for $T_M$.
The two other kernels are identical, except for the loop variable names, the \cppinline{vector_start}, and the leading dimension parameter for \cppinline{in} and \cppinline{out} passed to the BLAS call.
In contrast to the GPU implementation, the use of an optimized BLAS call is critical, as otherwise the kernel would not benefit from mixed-precision.
Due to the large cache size on CPUs, the vectors can be stored in the cache even for moderately large systems (e.~g. $100^3$), which means that other optimizations to ensure good cache locality are more important than the bandwidth.

\begin{listing}
\begin{minted}{cpp}
template<typename ValueType>
void omp_tensor_product_middle(
    gko::size_type num_rows, const ValueType* in, ValueType* out,
    const ValueType* B)
{
#pragma omp parallel for collapse(2)
  for (std::size_t k = 0; k < num_rows; ++k) {
    for (std::size_t i = 0; i < num_rows; ++i) {
      auto vector_start = k * num_rows * num_rows + i;
      cblas_gemv(CblasRowMajor, CblasNoTrans, num_rows, 
                 num_rows, ValueType(1.0), B, num_rows,
                 in + vector_start, num_rows, ValueType(0.0),
                 out + vector_start, num_rows);
    }
  }
}
\end{minted}
\caption{OpenMP kernel for applying $T_M$.}
\label{lst:cpu-tpm}
\end{listing}

\section{Numerical results}\label{sec:results}
\subsection{Hardware setup}

All benchmarks were run on the Noctua 2 cluster located at the Paderborn Center for Parallel Computing (PC2).
The GPU benchmarks were run on a single NVIDIA A100 40GB, and the CPU benchmarks were run on a two socket AMD EPYC Milan 7763 which contains 128 CPU cores in total.
The software stack uses GCC 12.2.0 as C++ compiler and CUDA 12.2.0.
Ginkgo is build from the commit \texttt{14ae29d} with the same GCC and CUDA versions.The BLAS backend used for the CPU runs is openBLAS 0.3.21.

\subsection{Test cases}

We consider two test problems in three dimensions, the heat equation as a prototype problem for a parabolic equation, and the advection equation, as a prototype for a hyperbolic equation. In both cases, we consider the domain $\Omega = [0,1]^3$. We employ second-order central finite difference stencils for the spatial discretization. 

\subsubsection{Heat equation}
We consider the 3D heat equation 
\[
\frac{\partial u}{\partial t} - \Delta u = g \text{ in } \Omega, \quad u = 0 \text{ on } \partial \Omega,
\]
with $g(\boldsymbol{x})=\sin(\pi x_1) \sin(\pi x_2) \sin(\pi x_3)$. The initial state ($t=0$) and analytical solution are given by
\[
u_a(\boldsymbol{x},t) = g(\boldsymbol{x})  \frac{ ( 1- e^{-3 \pi^2 t} )}{3 \pi^2}.
\]
When solving systems of the type as in \eqref{corelinear}, we use the preconditioner from \eqref{preconbasic} with the following settings
\begin{align*}
    A =  \mathbb{I}_n  +  \frac{\tau}{h^2} K_{1D}, \quad 
    B = C =  \frac{\tau}{h^2} K_{1D} 
\end{align*}
where $K_{1D} = \text{Tridiag}(-1,2,-1) \in \mathbb{R}^{n \times n} $ and $h=1/(n-1)$ is the spatial discretization. For the numerical tests of the heat equation, we exclude the 4s3pA method, as it exhibits unstable behavior for parameters of interest. This is expected, as outlined in the discussion on stability in~\Cref{sec:stability}.  
\subsubsection{Advection equation}
We consider the 3D Advection equation 
\[ 
\frac{\partial u}{\partial t} + \nabla \cdot \left(\bm{\beta}u\right) = 0 \text{ in } \Omega,
\]
with $\bm{\beta}=[1,1,1]^{\top}$ and periodic boundary conditions and the initial state given by 
\[
u(\bm{x}) = e^{ -100[(x_1-1/2)^2 + (x_2-1/2)^2 + (x_3-1/2)^2] }
\]
The corresponding 1D matrices used in constructing the preconditioner for systems of type \eqref{corelinear} are
\begin{align*}
    A =  \mathbb{I}_n  +  \frac{\tau}{2h} K_{1D} , \quad 
    B = C = \frac{\tau}{2h} K_{1D}  , \quad 
\end{align*}
where $K_{1D} = \text{Tridiag}(-1,0,1) + K^{(B)} \in \mathbb{R}^{n \times n} $ and $K^{(B)}_{1,n}=-1$, $K^{(B)}_{n,1}=1$ and $K^{(B)}_{ij}=0$ when $(i,j) \notin \{(1,n),(n,1)\}$. Here $K^{(B)}$ represents the 1-dimensional periodic boundary condition and we have $h=1/n$. 

\subsection{Linear solvers and setup of the numerical experiments}
For both test problems, we apply the preconditioners of the form \eqref{Pinvform}. This approach entails calculating the full spectral decomposition of the $n$-dimensional 1D matrices in setting up the preconditioner, the largest system considered is of dimension $N=n^3=64,000,000$, i.e., the preconditioner setup cost is the spectral decomposition of three matrices of at most dimension $400$---which is a negligible cost in relation to the solver. 

The combination of iterative linear solvers with mixed precision arithmetic deserves special attention, as solver convergence---and thus potential speedup---is affected by the introduction of low-precision computation. For this reason, we consider a range of different tolerances in stopping criteria in the CG and GMRES solvers. In all cases, we stop when either the absolute or the relative residual is less than the selected tolerance. In all tests we solve the equations to the final time $T_{end}=1/10$ and we vary the time steps $\tau \in \{ 10^{-1} , \  20^{-1} , \  40^{-1}, \  80^{-1} , \  160^{-1} , \  320^{-1} ,  \ 640^{-1} , \ 1280^{-1} \}$. {In all cases we set the maximum number of iterations in the GMRES/GC solvers to be 40.}

In this article, we focus on shared-memory performance. {The spatial discretization is constrained by the 40GB memory available on a single A100 GPU}. {In addition to the system matrix, the methods need to store the matrices $(a_{ij} \mathbb{I} + \tau K)$ and the Krylov basis for the GMRES solver.
For the mid-point rule, 4s3pA, and 4s3pB methods, only one additional matrix is necessary, as the diagonal entries of the respective Butcher tableaus $a_{ii}$ are identical.
However, the 4s3pC method requires four additional matrices to accompany the varying diagonal entries.
Thus, this method has a significantly higher memory requirement, which restricts the discretization size more than for the other methods.}
The spatial discretizations considered are $n \in \{200, \ 300, \ 400\}$, e.g., the spatial degrees of freedom (DOF) of the solution vector range from 8 to 64 million. These, in a 3D setting, relatively coarse spatial discretizations mean that the spatial error becomes dominating as $\tau$ decreases---thus this setting is not optimal for investigation of method error behavior with regard to $\tau$. This topic is the focus of other works, we refer the reader to \cite{PerformanceEvaluationRK}. This is however not an obstacle for the current objective of evaluating mixed precision performance for 3D PDE problems on different hardware. 

We use Ginkgo's built in profiler to separately time the following operations
\begin{itemize}
    \item Application of tensor operations from the left, middle and right in the form of $T_L$, $T_M$ and $T_R$ in \eqref{TensorDecomp}.
    \item Application of the full preconditioner in \eqref{Pinvform}.
    \item Application of the full preconditioned linear solver. We use CG in the case of the heat equation and GMRES for the advection equation.
\end{itemize}

\subsection{Combined error plots}
Before focusing on the behavior of mixed-precision implementations in Ginkgo, we validate the used methods by including the behavior of the discretization errors. As discussed in the previous section, the current setting of a single GPU with second order spatial discretization means that we expect the spatial error to be the dominating component as $\tau$ decreases. 

The results for both equations are presented in \Cref{fig:CombinedErrors}.
Only the cases for solver tolerance tol=$10^{-3}$ is considered here, as this already captures the behavior sufficiently as the error curves look essentially the same for all tolerances. 
For the heat equation in Figure~\ref{fig:CombinedErrors:Heat}, we do not see a significant difference between single and double precision for the midpoint rule with correction (M), the error behaves as $\mathcal{O}(\tau^2)$ as expected. The 4s3pB method (B) follows the third order convergence and the single/double results diverge as $\tau$ and the error decreases. Here the single precision results in a lower error, but as both precisions result in errors $\sim 10^{-7}, \  10^{-8}$, the difference can be explained by rounding errors resulting in a solution closer to the analytical values. The 4s3pC method (C) displays some order reduction, as the observed error behavior is in-between second and third order, with minor differences due to precision. 

\begin{figure}[tbp]
\centering
{\includegraphics[width=0.6\textwidth]{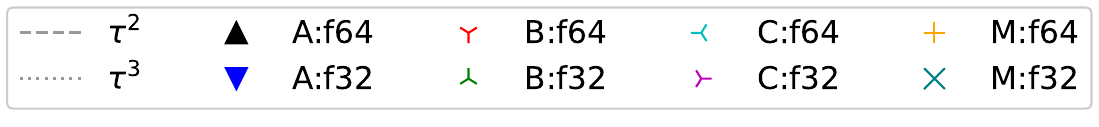}} \\
\subfigure[Heat eq. solver tolerance $10^{-3}$]{\includegraphics[width=0.49\textwidth]{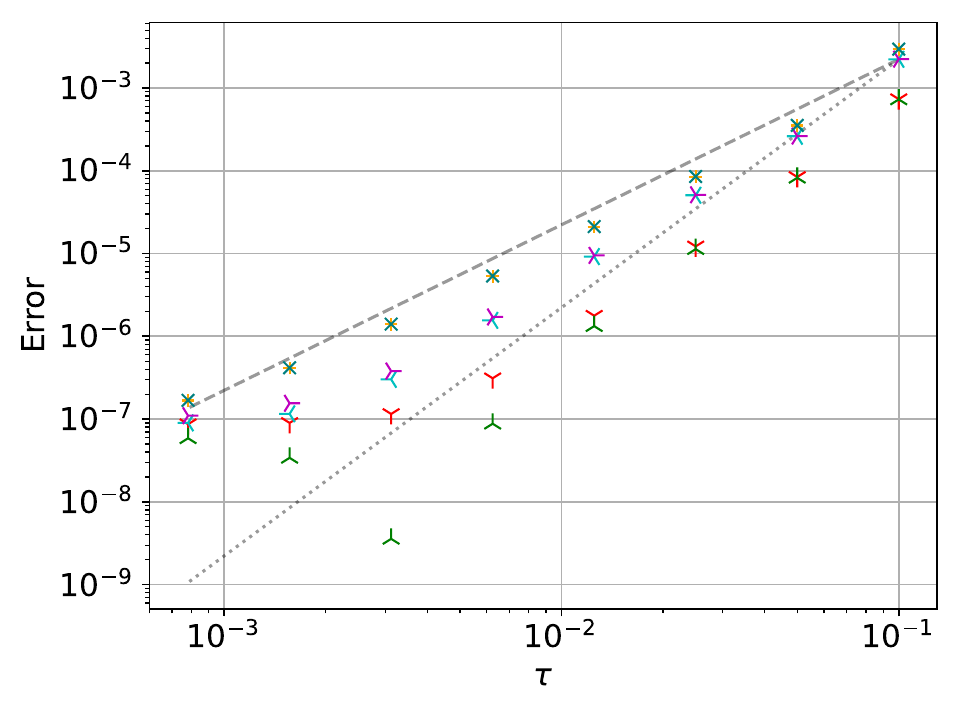}\label{fig:CombinedErrors:Heat}}
\subfigure[Advection eq. solver tolerance $10^{-3}$ ]{\includegraphics[width=0.49\textwidth]{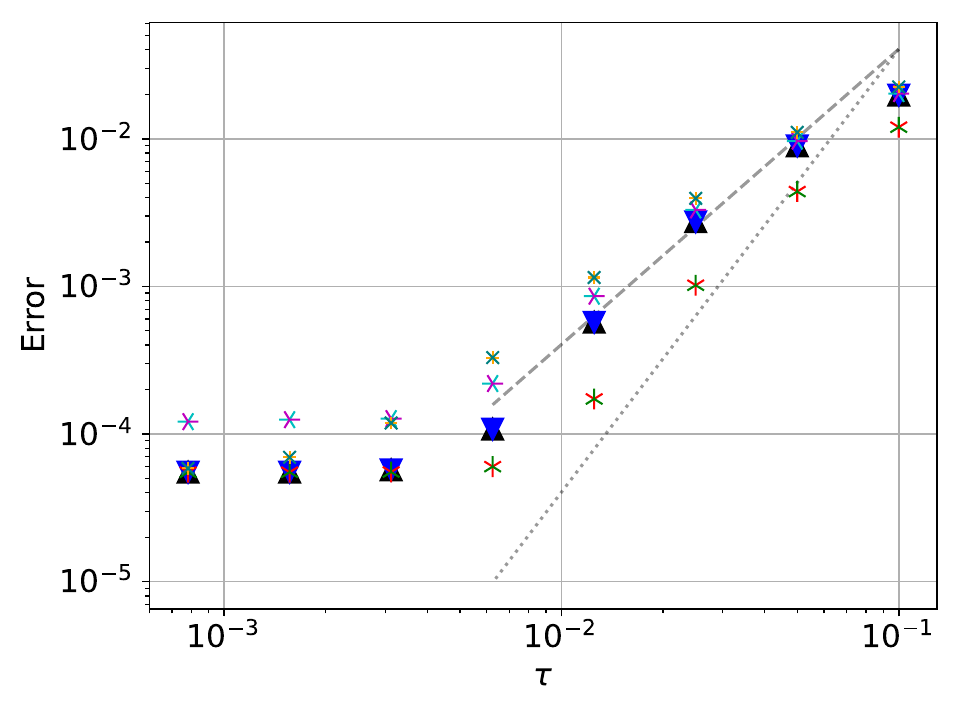}\label{fig:CombinedErrors:Advection}}
\caption{{Error at the final time for Methods A, B, C and midpoint rule with correction.}}
\label{fig:CombinedErrors}
\end{figure}

The errors for the advection equation are presented in \Cref{fig:CombinedErrors:Advection}. Compared to the heat equation we see that the final error of the methods plateaus at significantly higher values. The observed errors are consistent with accumulated spatial error components, $h^2=(1/300)^2$ for methods A, B and M and $h^2=(1/200)^2$ for method C. The order reduction is more pronounced, with methods A and B displaying an order between two and three, while methods C and M are closer to second order. As the error plateaus are reached at fairly large time steps, it is possible that the methods are not yet operating in their asymptotic regime. As in the case of the heat equation, the tolerance used in the iterative solver for the implicit steps seems to have negligible effect when combined with the explicit corrector steps. 

An important conclusion from this section is the fact that reducing the precision in certain components of the mixed-precision Runge--Kutta time stepping methods does not reduce the accuracy of the calculations for the chosen parameters. In the next sections, we will consider the speedup obtained by mixing singe and double precision which is, hence, obtained without loss in accuracy.

\subsection{Speedup of tensor components}\label{Sec:Tensor_comp}

\Cref{fig:CombinedTensorspeedup} shows the speedup of the three tensor-product applications when comparing single to double precision, which form the preconditioner, when running on GPUs. The run times are collected by solving the problem until the final time 0.1 with  the respective time integrators and normalizing the run time over the total number of solver iterations.
The average over the methods A, B, C, and mid-point, is shown, as well as the speedup for each method individually in low opacity---illustrating the spread. 
Overall, these components achieve speedups between 1.25 and 2.5.
The individual behavior of the three kernels, however, is not consistent for both equations.

\begin{figure}[tbp]
\centering
{\includegraphics[width=0.35\textwidth]{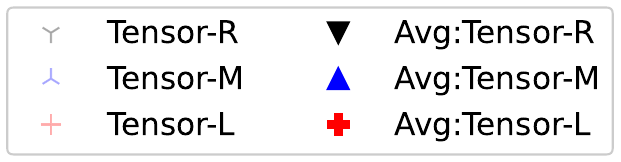}} \\
\subfigure[GPU: heat equation (tol=$10^{-3}$)]{\includegraphics[width=0.49\textwidth]{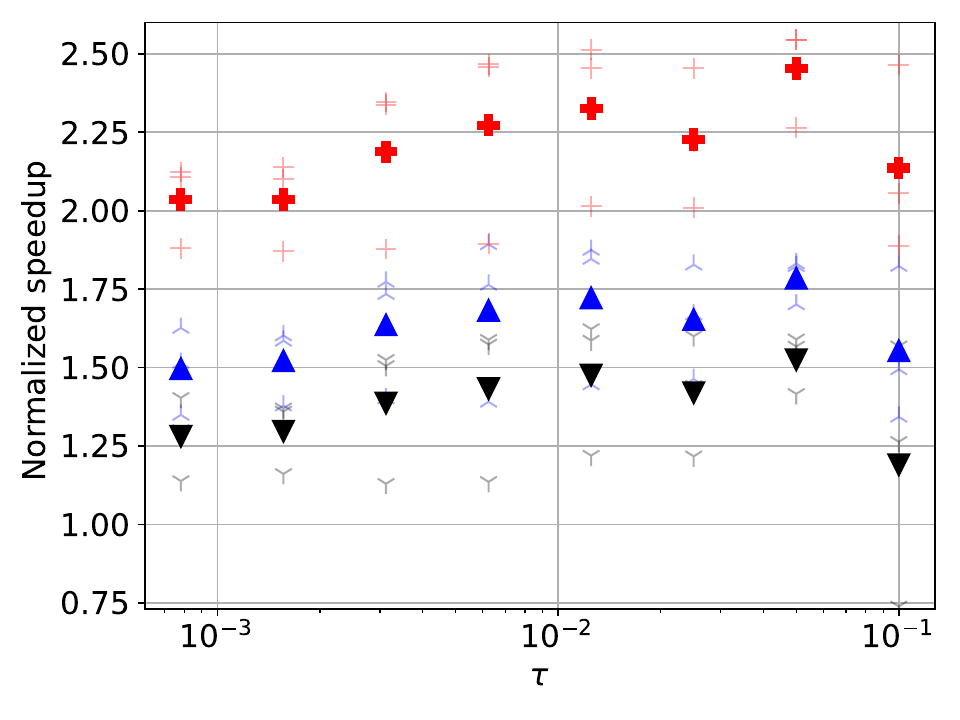}}
\subfigure[GPU: Advection equation (tol=$10^{-3}$)]{\includegraphics[width=0.49\textwidth]{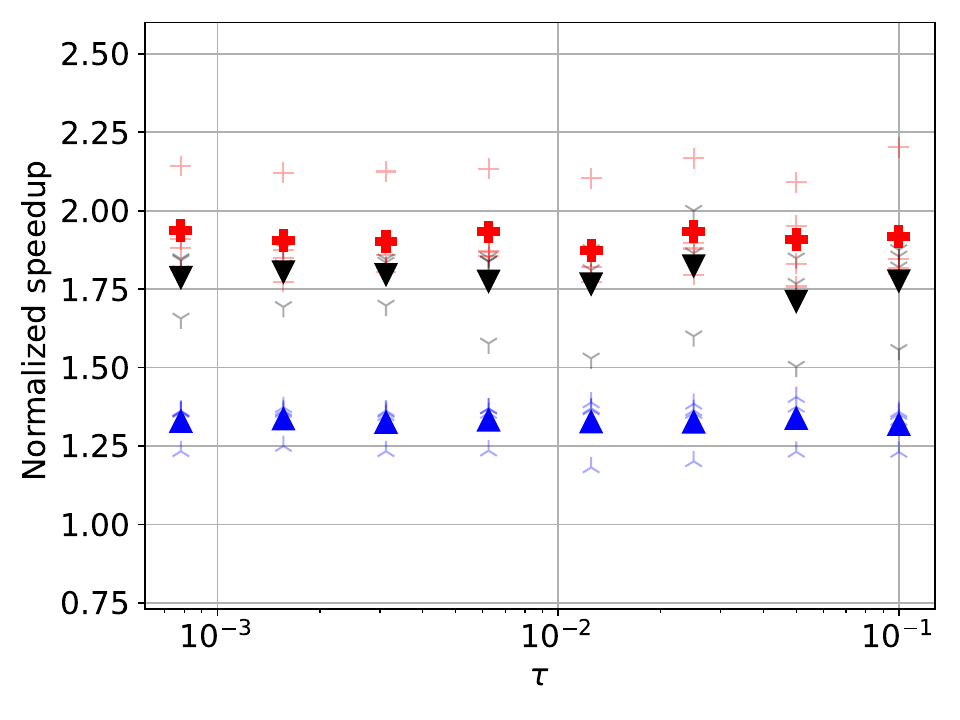}}
\caption{GPU: Normalized (over iterations) speedup and averages for tensor operations $T_{R,M,L}$.}
\label{fig:CombinedTensorspeedup}
\end{figure}

First, the Tensor-L kernel shows the best speedups for both equations.
It is a purely main memory bandwidth limited kernel, and as such benefits the most from the reduced data volume.
In the heat equation case, it even shows more than ideal speedup.
One reason could be the very few number of iterations for double precision, see \Cref{tbl:iterations}, as they could introduce higher variance and leave the GPU unutilized for too long.
Thus, the advection case should be considered as the fairer comparison.
It should also be noted, that this kernel is overall the slowest one.
Compared to Tensor-R, and Tensor-M it is slower by a factor of $2\times$ and $4\times$, respectively, as can be seen from \Cref{fig:solvertime}.

The other two kernels show less speedup, although it differs by the equation to solve.
The Tensor-R kernels achieves close to $1.5\times$ speedup for the heat equation, and $2\times$ speedup for the advection equation.
The Tensor-M kernel shows the opposite behavior.
It achieves $1.5$--$1.8\times$ speedup for the heat equation, but less than $1.5\times$ speedup for the advection equation.
Compared to the Tensor-L kernel, both kernels are more compute bound.
The access pattern with unit- and $n$-stride for the Tensor-R and Tensor-M kernel, respectively, leads to efficient caching, as a single thread block shares only a few 1D blocks of the input vectors.
Thus, the issued warps can easily reuse data from the L2 cache.
The differences between the heat and advection equations for these two kernels cannot be explained with the measured data.
As they are more compute bound, the memory-focused optimizations and analysis proved insufficient.
Further studies could investigate optimizations for compute bound codes, such as vectorizing the innermost loop, or targeting the tensor cores of the GPUs.
The study~\cite{cui-tensor-products} shows promising optimizations, especially within the context of lower precision and tensor-product applications.

The speedups of the tensor kernels running on a CPU are presented in~\Cref{fig:CombinedTensorspeedupCPU}.
For the heat equation, the behavior is similar to the GPU case; the kernels achieve good speedup through using single precision over double precision, and the best speedup, close to $2\times$ is achieved for the Tensor-L kernel.
The other two kernels achieve less speedup, at most $1.5\times$.
However, the speedup in the advection equation case is considerably worse than on GPUs.
As the complex-valued computations are mostly handed off to the BLAS library, the inefficiency might be caused from fewer optimizations on the particular BLAS calls used here, i.e., \cppinline{gemv} with complex values and large strides.
Additionally, the CPU implementation might suffer from inefficient use of the non-uniform memory access (NUMA) domains, since the full node is used, as mentioned in~\cite[Chapter~8]{hager_introduction_2017}.

\begin{figure}[tbp]
\centering
{\includegraphics[width=0.35\textwidth]{plots/legend_GPU_Heat_combined_tensor_speedups_1em3.pdf}} \\
\subfigure[CPU: heat equation (tol=$10^{-3}$)]{\includegraphics[width=0.49\textwidth]{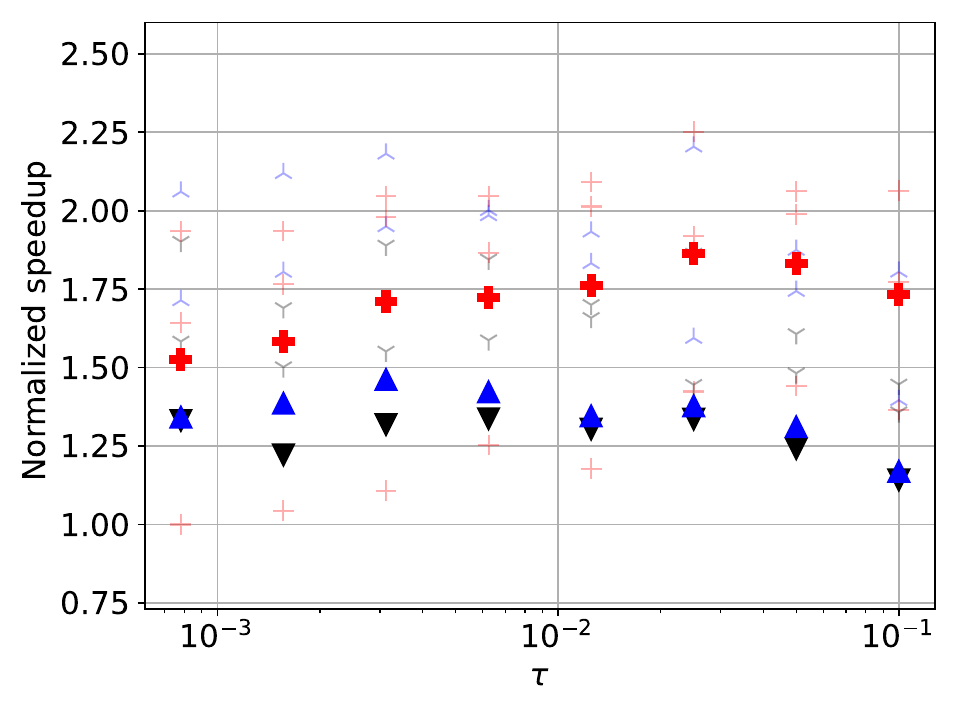}}
\subfigure[CPU: Advection equation (tol=$10^{-3}$)]{\includegraphics[width=0.49\textwidth]{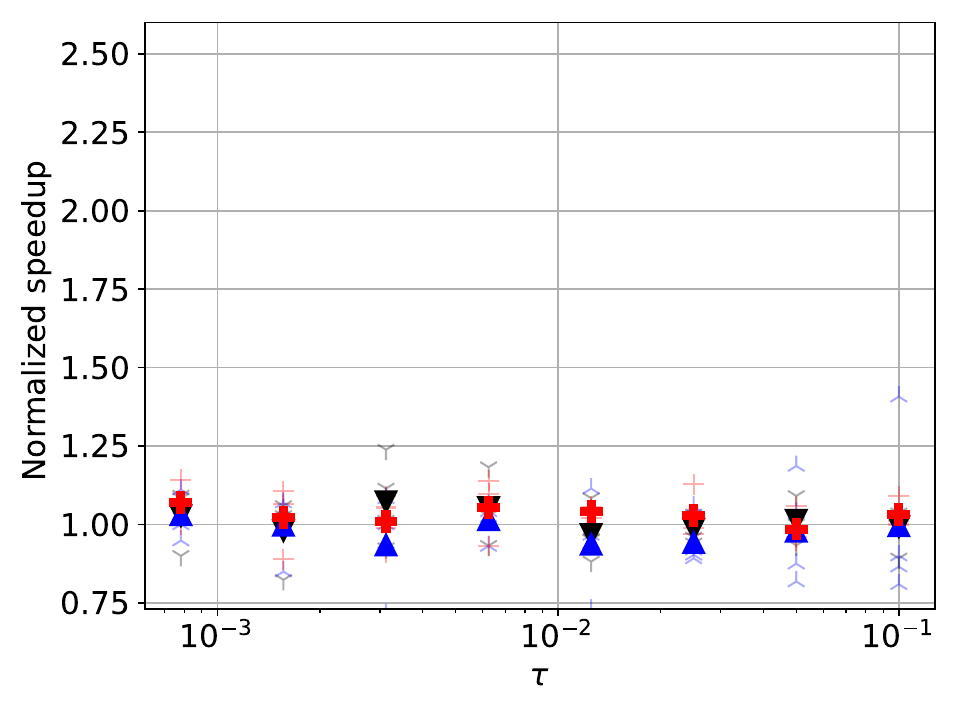}}
\caption{CPU: Normalized (over iterations) speedup and averages for tensor operations $T_{R,M,L}$.}
\label{fig:CombinedTensorspeedupCPU}
\end{figure}

\subsection{{Speedup of linear solver and preconditioner}}\label{Sec:speedupsolver}

{The speedup over the full solver application is not only given by the speedup of the individual tensor product applications.
The lower precision can also affect the solver convergence rate.
To illustrate, \Cref{tbl:iterations} shows the number of solver iterations until convergence averaged over all time steps for the Method C with $\tau = 0.0125$.
While the chosen preconditioner is ideal for the heat equation solved in double precision leading to convergence in one iteration, the reduced precision introduces additional rounding errors, which lead to a deterioration of the convergence. Thus, the single precision GPU solver for the heat equation requires at least twice as many iterations.
This effect is reduced for the CPU solver, most likely due to the use of optimized BLAS routines, which can also reduce rounding errors as described in~\cite{Higham_Mary_2022}. In all cases, however, the number of iterations increase starting at tolerances about $10^{-5}$, and the solver convergence becomes very slow or is even lost as the tolerance approaches single precision machine epsilon.
For the advection equation, the preconditioner is less ideal, leading to smaller differences in the iteration count.}

\begin{table}[tpb]
    \small
    \centering
    \begin{tabularx}{\textwidth}{lrrrrrrrr}\toprule
& \multicolumn{4}{c}{CPU} & \multicolumn{4}{c}{GPU} \\\cmidrule(lr){2-5}\cmidrule(lr){6-9}
& \multicolumn{2}{c}{Advection} & \multicolumn{2}{c}{Heat} & \multicolumn{2}{c}{Advection} & \multicolumn{2}{c}{Heat} \\\cmidrule(lr){2-3}\cmidrule(lr){4-5}\cmidrule(lr){6-7}\cmidrule(lr){8-9}
Tolerance & \multicolumn{1}{c}{Float32} & \multicolumn{1}{c}{Float64} & \multicolumn{1}{c}{Float32} & \multicolumn{1}{c}{Float64} & \multicolumn{1}{c}{Float32} & \multicolumn{1}{c}{Float64} & \multicolumn{1}{c}{Float32} & \multicolumn{1}{c}{Float64} \\\hline
0.001 & 2.5 & 2.5 & 1.0 & 1.0 & 2.5 & 2.5 & 1.61 & 1.14\\
0.0001 & 3.0 & 3.0 & 1.46 & 1.0 & 2.5 & 3.0 & 1.64 & 1.0\\
1.0e-5 & 4.09 & 4.0 & 2.14 & 1.0 & 4.0 & 4.0 & 2.0 & 1.0\\
1.0e-6 & 4.99 & 4.5 & 2.93 & 1.0 & 4.5 & 4.5 & 2.67 & 1.0\\
1.0e-8 & --- & 6.5 & 5.75 & 1.0 & --- & 6.5 & 5.75 & 1.0\\\bottomrule
\end{tabularx}
    \caption{The number of solver iterations until convergence for the Method C. The value is averaged over all time steps with $\tau = {1/640}$. Note that the GPU advection case uses $200^3$ DOFs, while all other cases use $300^3$ DOFs. (---) denotes failure to converge.}
    \label{tbl:iterations}
\end{table}

Before studying the overall solver performance, we also consider the share that the individual operations of the solver have on the total run time. These are shown in Figure \ref{fig:solvertime}, we only display one case, namely Method C with a solver tolerance of $10^{-3}$. We illustrate the time distribution of Method C (cf.~\Cref{tbl:iterations}), as the relative timings are nearly identical across solver tolerances and even heat/advection, we include only the double precision heat example with the tolerance $10^{-3}$. As can be seen from the timing plots, the tensor operation Tensor-L, is dominating the total compute time on GPU, while in the CPU case the times for the three operations are far more equal. In the CPU case, the other operations in the solver also take up far more of the total run time when compared against the GPU case.

\begin{figure}[tbp]
\centering
\subfigure[GPU: Heat equation]{\includegraphics[width=0.49\textwidth]{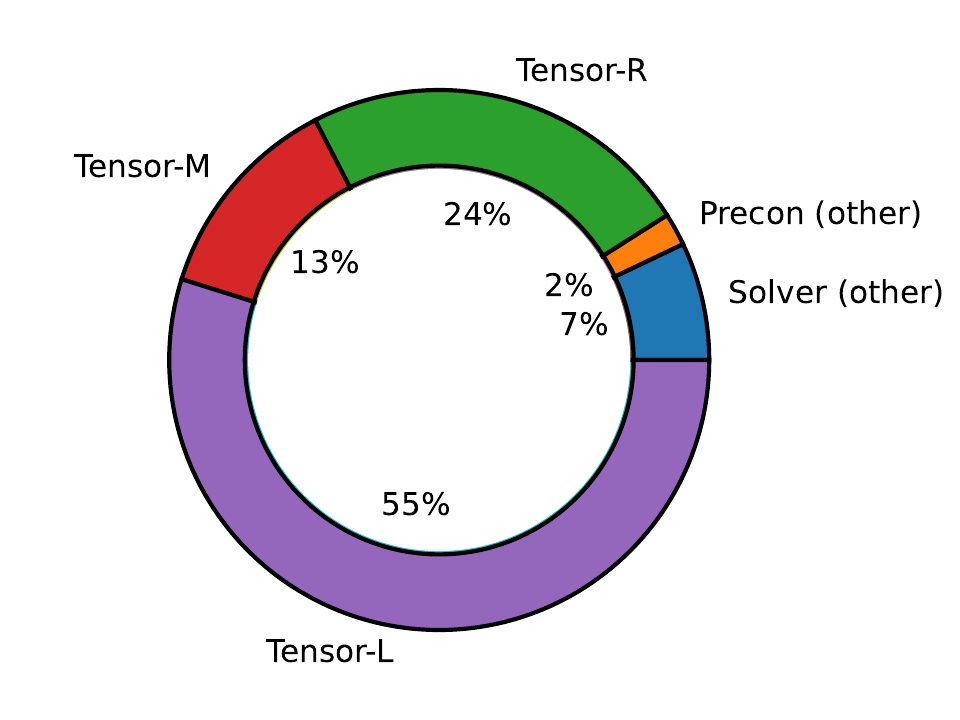}}
\subfigure[CPU: Heat equation]{\includegraphics[width=0.49\textwidth]{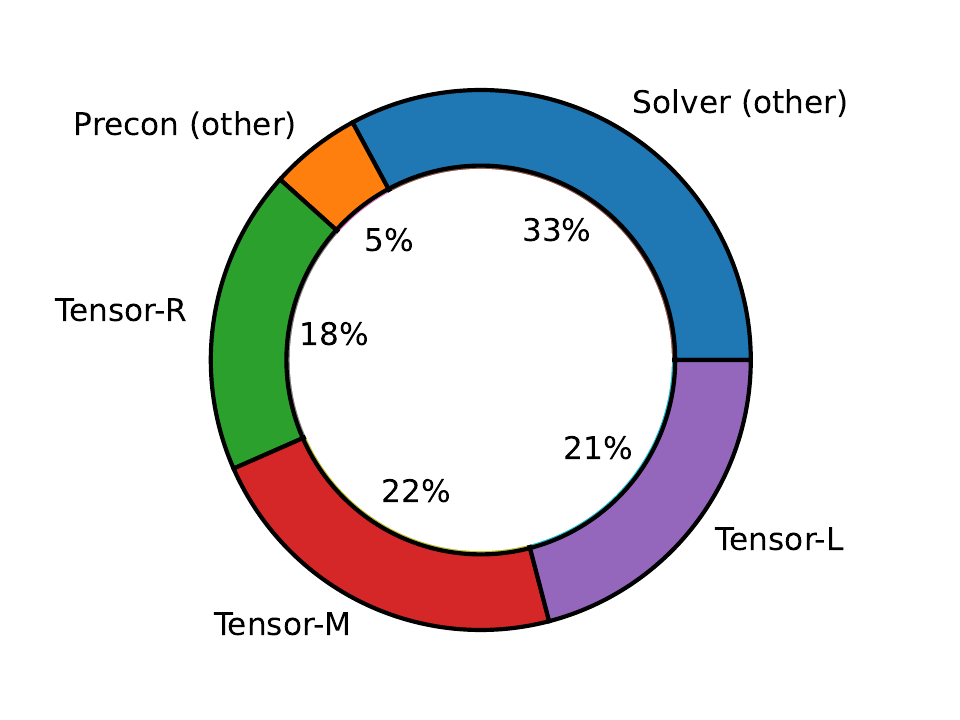}}
\caption{Proportion of Time spent in different parts of the solver. (Method C, double precision tol=$10^{-3}$, timings across all time-steps).}
\label{fig:solvertime}
\end{figure}

\Cref{fig:GPU_CombinedSpeedupHeat} and \Cref{fig:GPU_CombinedSpeedupConvection} present the speedup of running the linear solver in single precision over double precision on GPU, and \Cref{fig:CPU_CombinedSpeedupHeat} and \Cref{fig:CPU_CombinedSpeedupConvection} on CPU for the heat and advection case, respectively.
The figures represent the gain over the solver application across all time steps.
This means they are not normalized over the number of iterations, and instead, the speedups already consider the differences in convergence that might occur due to the different precisions.
Note, that we focus on the solver and time stepping method while parts related to the spatial discretization, e.g., the assembly of the right-hand-side, are not considered in this study. These components are not crucial to the mixed precision behavior and in general will depend on the discretization used.

\subsubsection{Heat}

For the GPU case, the speedup of using single precision for the heat equation shows a clear upwards trend in \Cref{fig:GPU_CombinedSpeedupHeat}.
For larger tolerances, the speedup increases with decreasing time step size, and reaches a speedup of up-to $1.4\times$ for the smallest time step.
However, the smallest tolerance $10^{-8}$ shows a slowdown\footnote{This tolerance is already at the machine epsilon for single precision.}.
Both behaviors are linked to the difference in convergence of the single and double precision solvers, as mentioned above.
This becomes less of an issue as the time-step size decreases, {as, for a fixed $n$, the condition number of the preconditioned system goes to one in the limit $\tau \rightarrow 0$.}
If the required iteration numbers approach each other, the gain for the individual preconditioner components, as shown in \Cref{fig:CombinedTensorspeedup}, become more pronounced, leading to better speedups.

The CPU results are shown in \Cref{fig:CPU_CombinedSpeedupHeat}, where the general trend is similar to that of the GPU case. Speedups tend to improve as the time step decreases and the choice of small solver tolerances is detrimental to speedup of the low precision implementation due to the increase in iterations observed when solving to tolerances close to corresponding machine precision.

\begin{figure}[tbp]
\centering
{\includegraphics[width=0.7\textwidth]{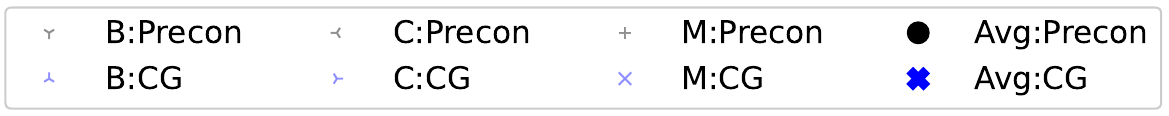}} \\
\subfigure[tol=$10^{-3}$]{\includegraphics[width=0.49\textwidth]{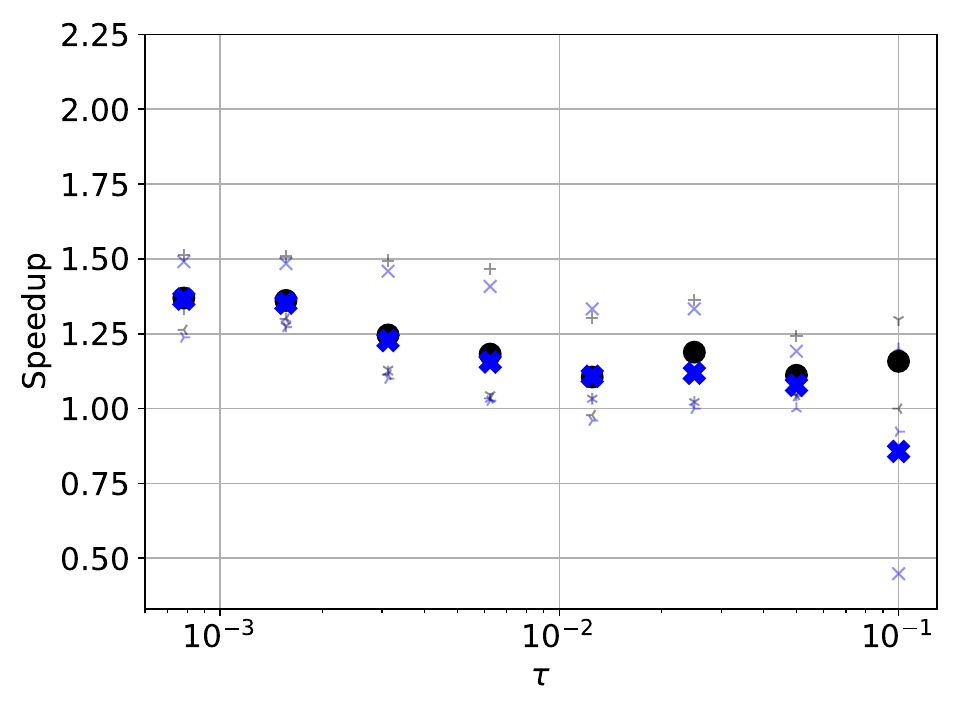}}
\subfigure[tol=$10^{-4}$]{\includegraphics[width=0.49\textwidth]{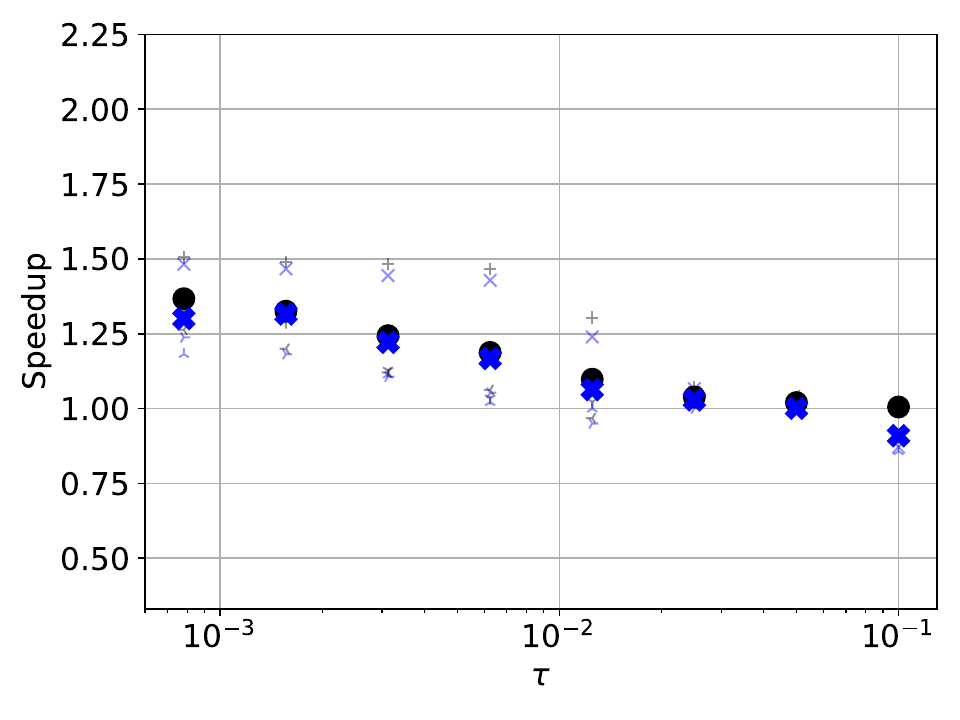}}
\subfigure[tol=$10^{-5}$]{\includegraphics[width=0.49\textwidth]{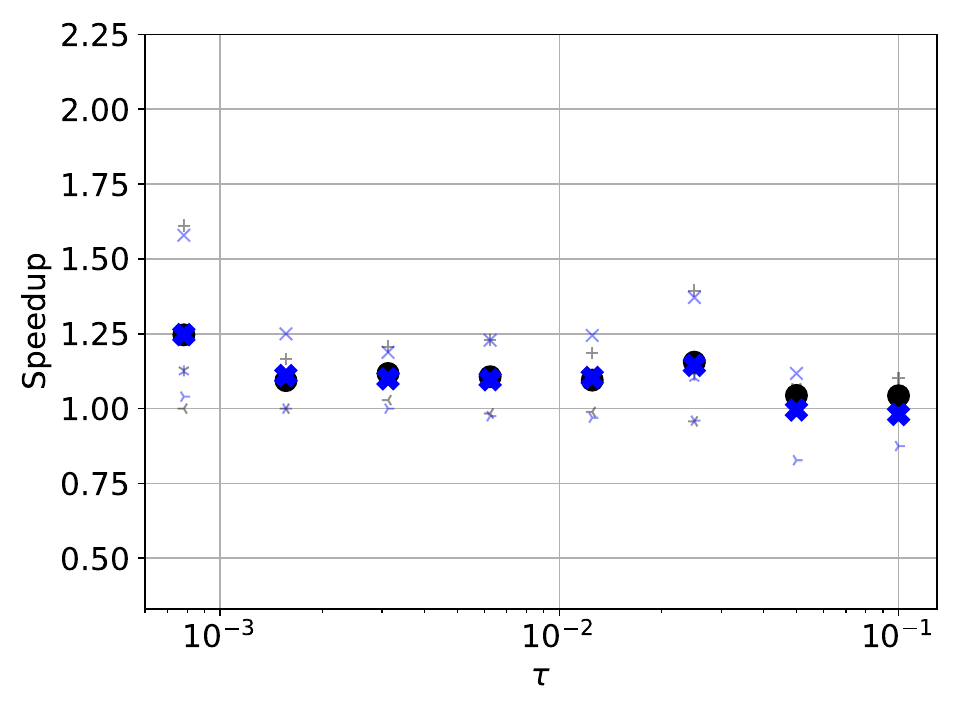}}
\subfigure[tol=$10^{-8}$]{\includegraphics[width=0.49\textwidth]{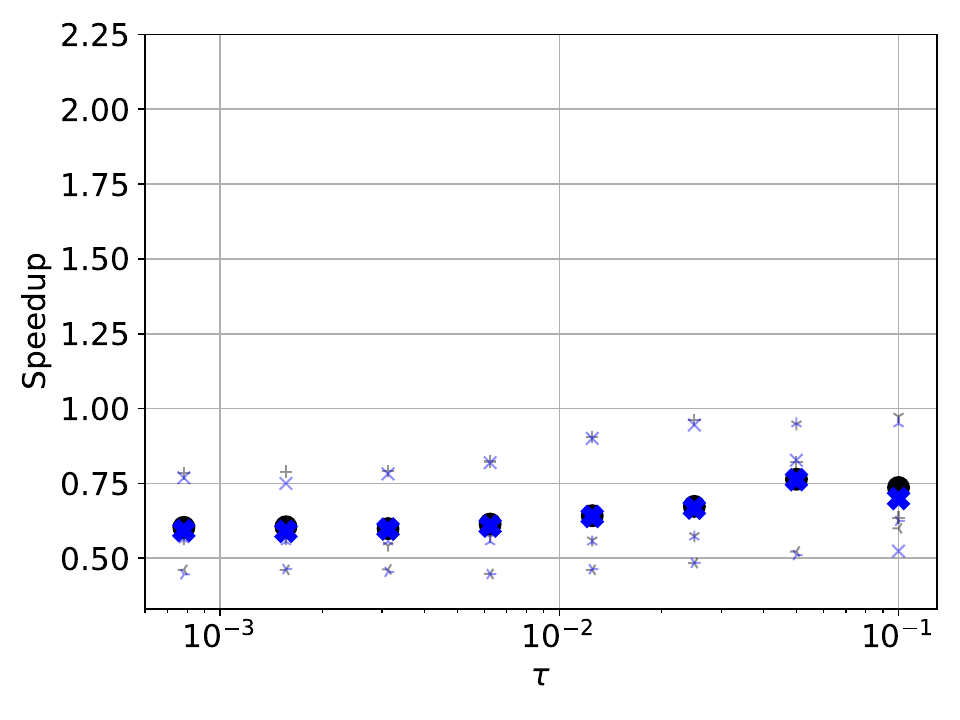}}
\caption{GPU: Heat equation: combined speedup and averages. Methods B, C and midpoint rule with correction.}
\label{fig:GPU_CombinedSpeedupHeat}
\end{figure}

\begin{figure}[tbp]
\centering
{\includegraphics[width=0.7\textwidth]{plots/legend_GPU_Heat_1em3.pdf}} \\
\subfigure[tol=$10^{-3}$]{\includegraphics[width=0.49\textwidth]{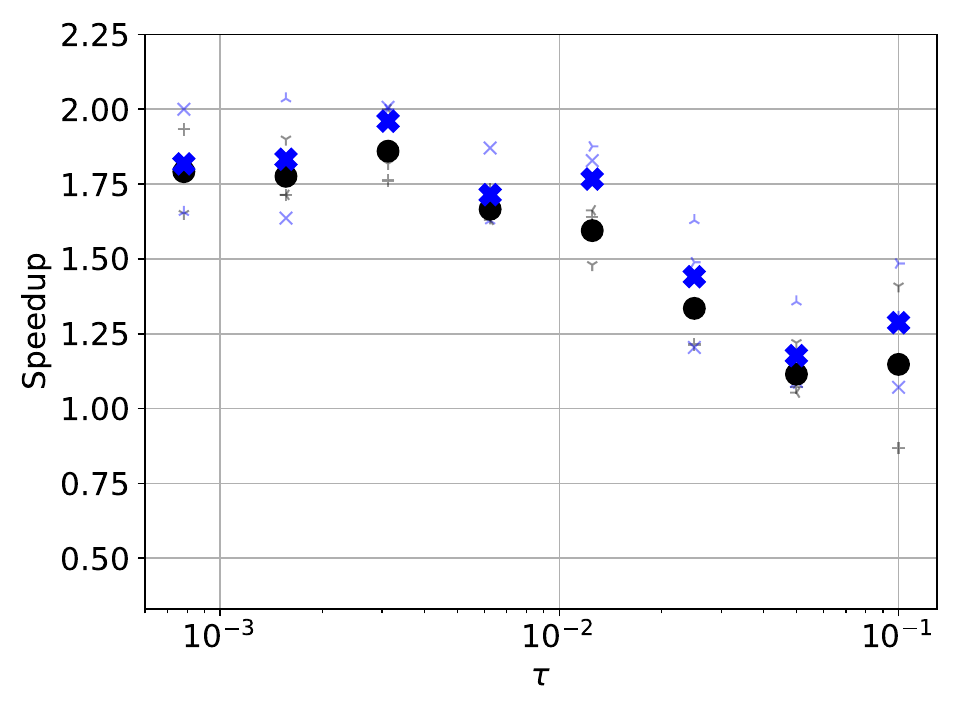}}
\subfigure[tol=$10^{-4}$]{\includegraphics[width=0.49\textwidth]{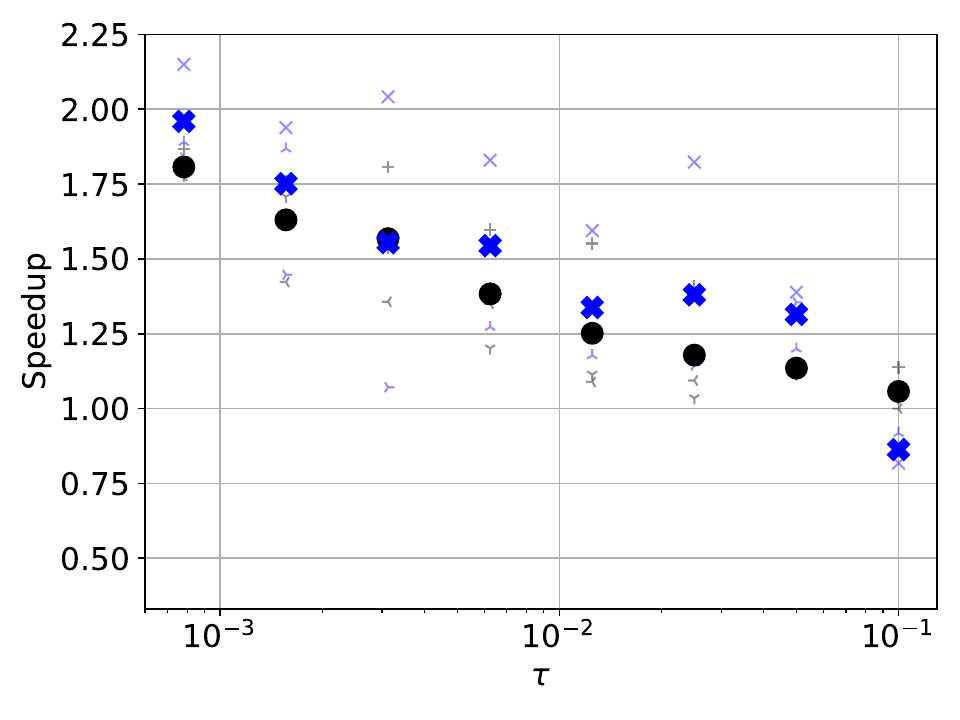}}
\subfigure[tol=$10^{-5}$]{\includegraphics[width=0.49\textwidth]{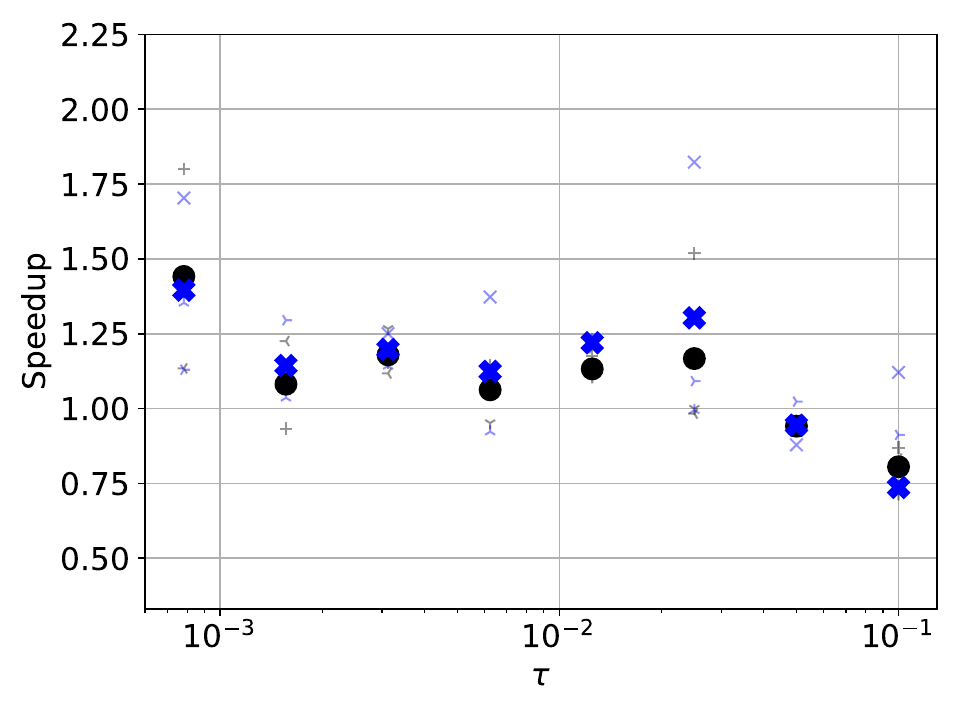}}
\subfigure[tol=$10^{-8}$]{\includegraphics[width=0.49\textwidth]{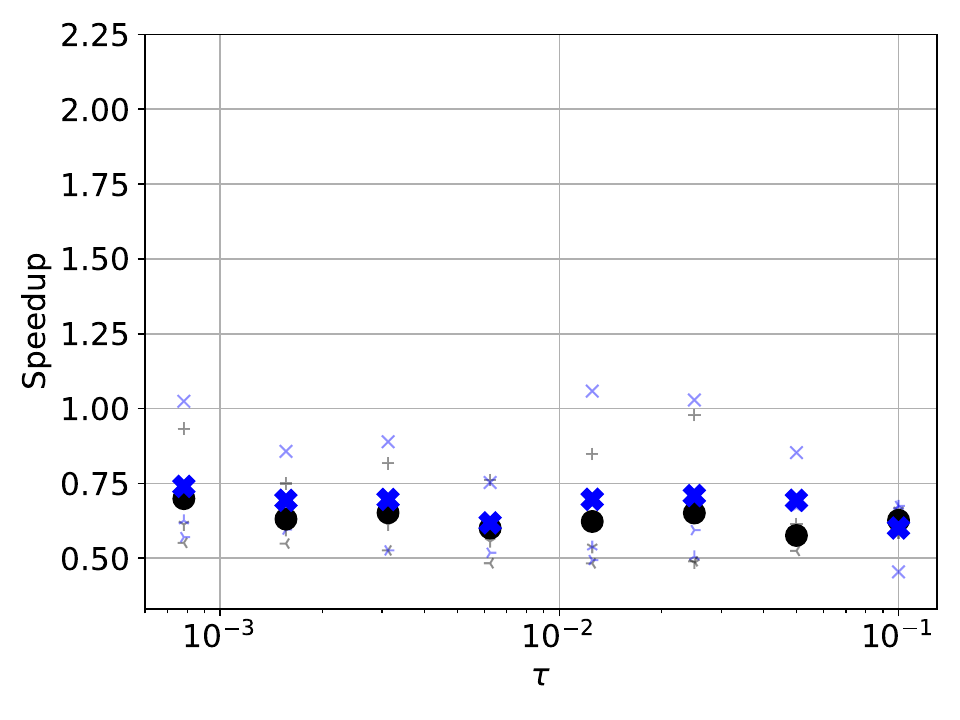}}
\caption{CPU: Heat equation: combined speedup and averages. Methods B, C and midpoint rule with correction.}
\label{fig:CPU_CombinedSpeedupHeat}
\end{figure}

\subsubsection{Advection}

In contrast to the heat equation, the GPU speedups for the advection equation shown in \Cref{fig:GPU_CombinedSpeedupConvection} are more consistent. For the considered tolerances, the speedup is mostly between $1.8-2.0\times$, with some outliers around $1.4-1.6\times$.
As shown in \Cref{tbl:iterations} both the single and double precision solvers converge at the same rate, which leads to a speedup independent of the time step size. Additionally, better speedups compared to the heat equation are possible due to higher memory requirements. The advection case requires using the GMRES solver, which stores additional Krylov basis vectors, as well as complex numbers, which directly doubles the memory requirements compared to the heat equation. 

The CPU case is illustrated in \Cref{fig:CPU_CombinedSpeedupConvection}. As foreshadowed in the results outlined in Section \ref{Sec:Tensor_comp}, at best we achieve only modest speedups of about $1.2\times$, with several methods resulting in a slower wall-clock times for the low precision calculations. In contrast to all other test cases, cases of less than $1$ speedups is observed for all considered tolerances.

\begin{figure}[tbp]
\centering
{\includegraphics[width=0.875\textwidth]{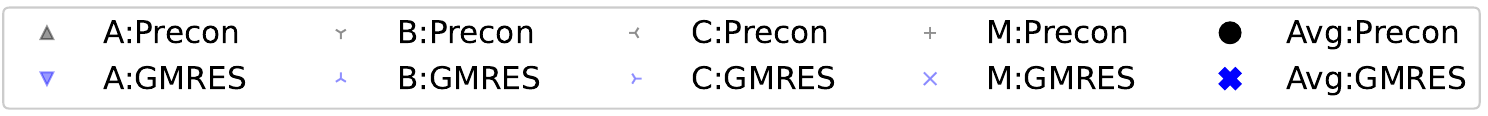}} \\
\subfigure[tol=$10^{-3}$]{\includegraphics[width=0.49\textwidth]{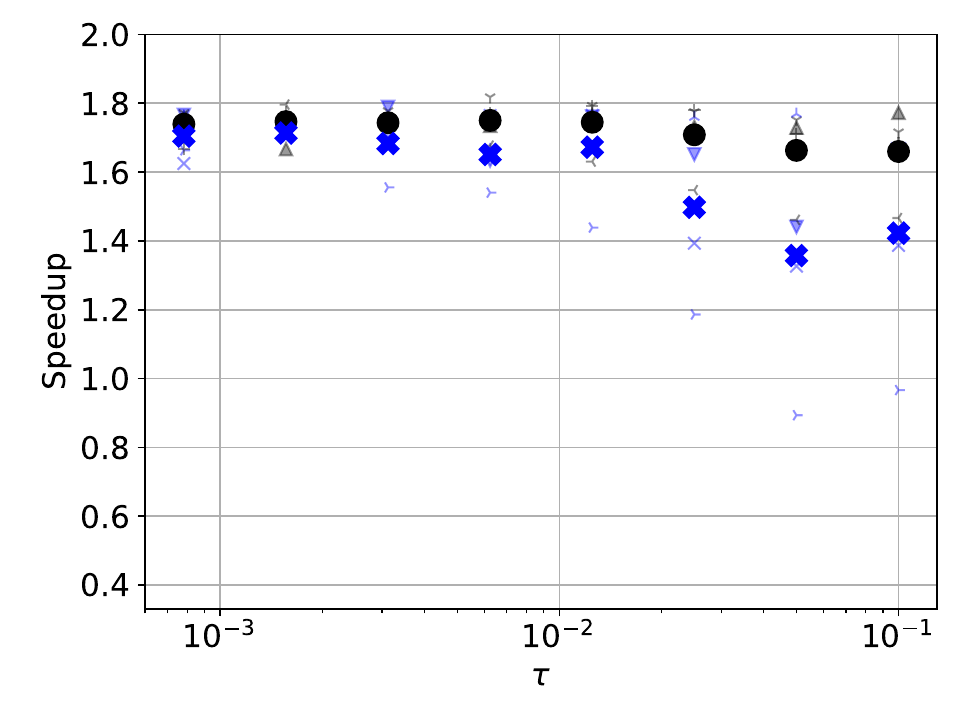}}
\subfigure[tol=$10^{-4}$]{\includegraphics[width=0.49\textwidth]{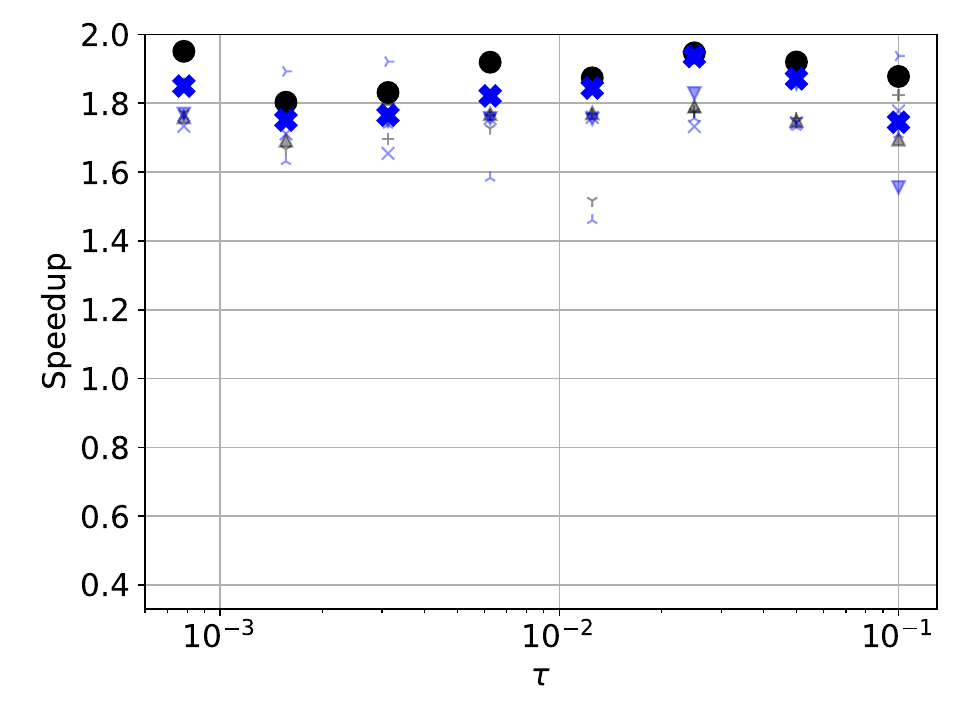}}
\subfigure[tol=$10^{-5}$]{\includegraphics[width=0.49\textwidth]{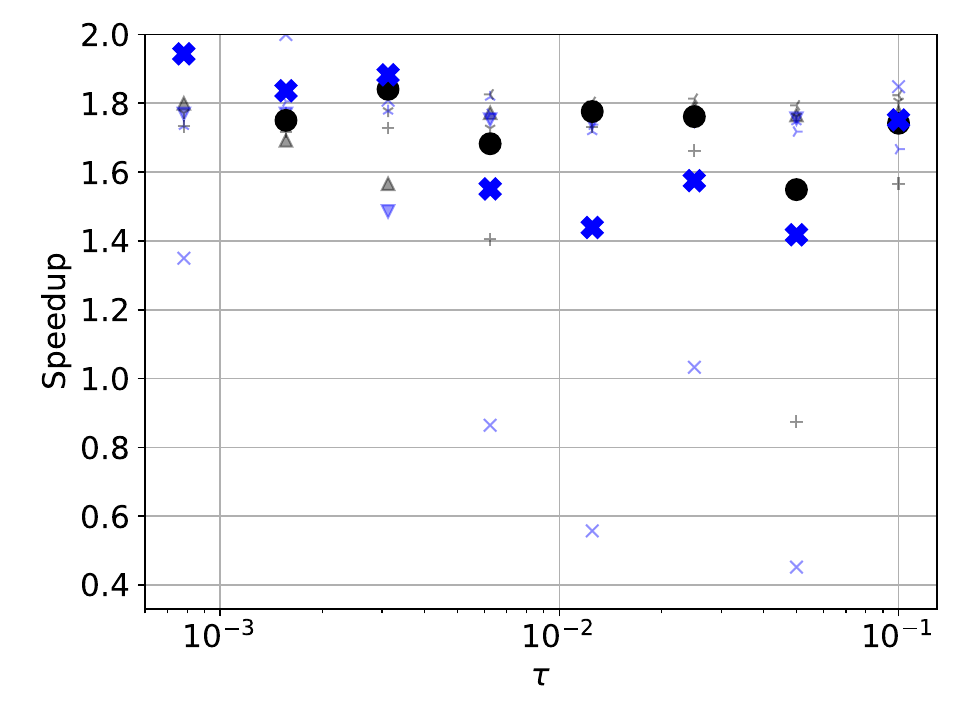}}
\subfigure[tol=$10^{-6}$]{\includegraphics[width=0.49\textwidth]{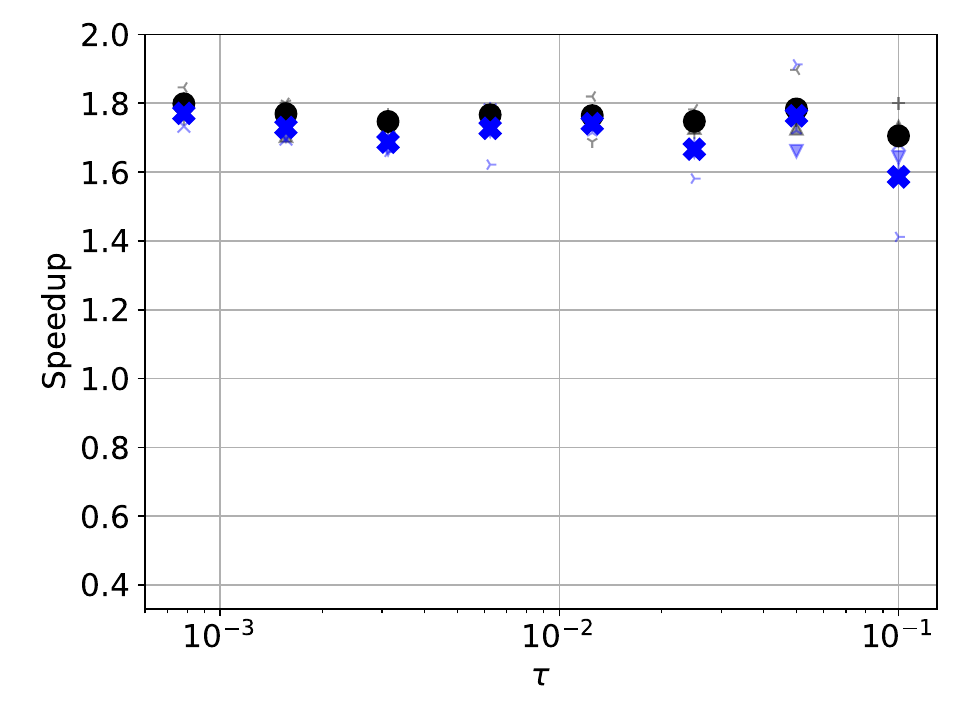}}
\caption{GPU: Advection equation: combined speedup and averages. Methods A, B, C and midpoint rule with correction.}
\label{fig:GPU_CombinedSpeedupConvection}
\end{figure}

\begin{figure}[tbp]
\centering
{\includegraphics[width=0.875\textwidth]{plots/legend_GPU_Advec_1em3.pdf}} \\
\subfigure[tol=$10^{-3}$]{\includegraphics[width=0.49\textwidth]{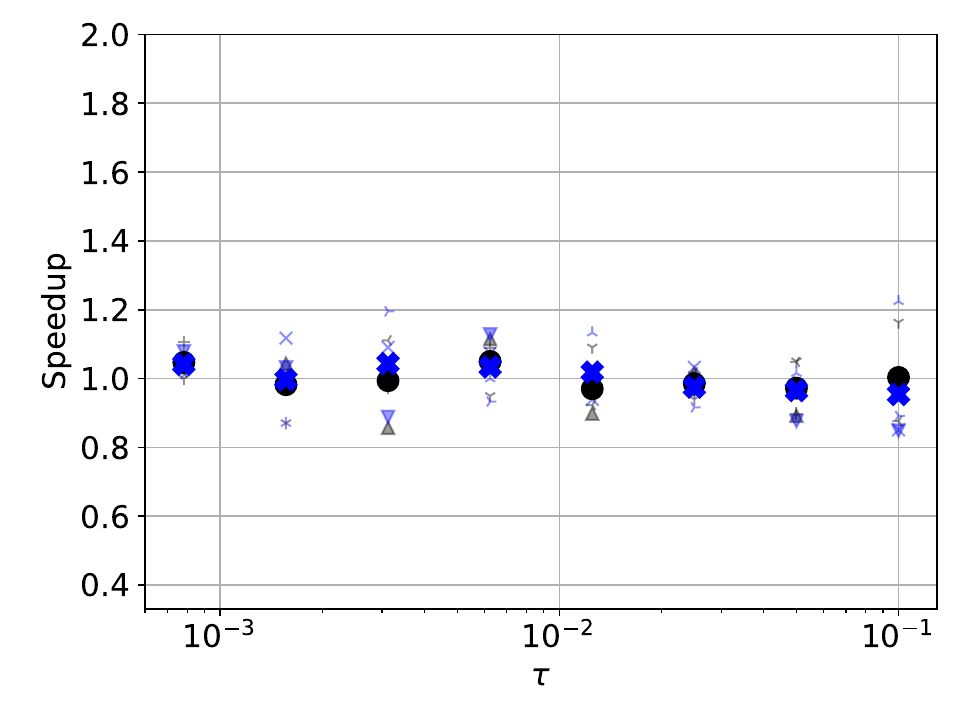}}
\subfigure[tol=$10^{-4}$]{\includegraphics[width=0.49\textwidth]{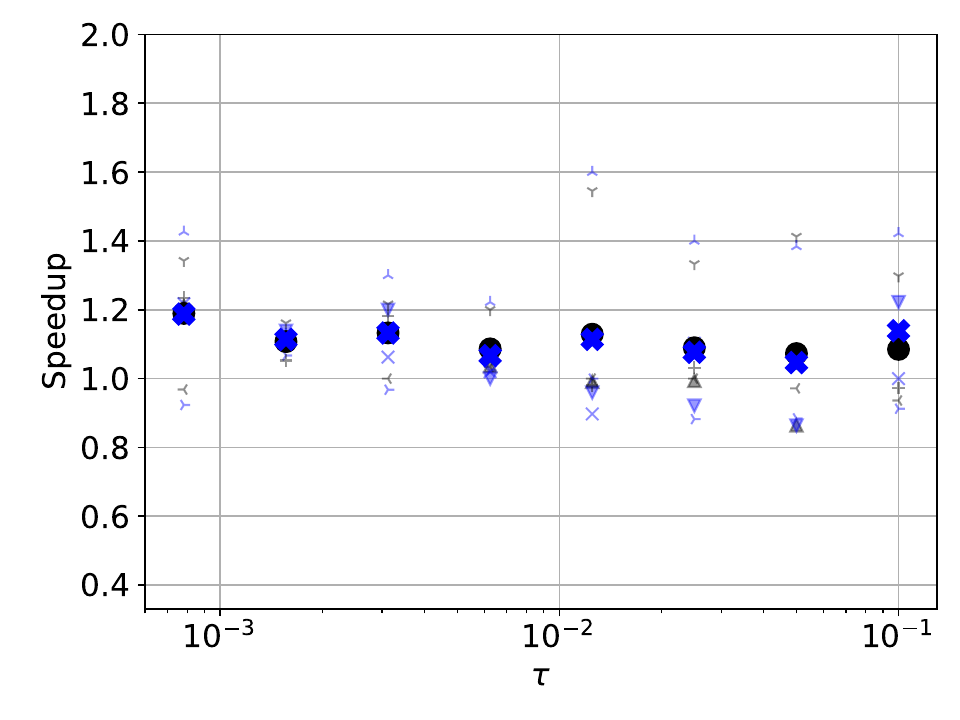}}
\caption{CPU: Advection equation: combined speedup and averages. Methods A, B, C and midpoint rule with correction.}
\label{fig:CPU_CombinedSpeedupConvection}
\end{figure}

\section{Summary and conclusions}\label{sec:conclusions}

In this paper, we investigated the performance of mixed-precision for solving time-dependent PDEs with implicit and partially implicit Runge--Kutta methods. Based on the previous work by Grant~\cite{Grant2022}, the explicit parts of the time-stepping methods were treated in high precision and low precision was used to solve the linear systems for the implicit parts.
The solvers were implemented using the Ginkgo library, which enables straightforward switching between precisions and hardware, and were run on both modern GPUs and CPUs. 
Two PDEs with different characteristics were considered, namely the heat and advection equation.

Our numerical study shows that the memory-bound building blocks reach optimal speedup---or even more---while the speedup is reduced for operations that are compute bound. The speedup obtained by the total solver can be reduced by an increase in iterations necessary for convergence. The study shows that the choice of the solver tolerance is crucial and needs to be chosen with care when using reduced precision. Also care has to be taken when implementing matrix-vector products in lower precision to reduce roundoff errors. We observe that the BLAS routine \texttt{gemv} achieves this very well in the real case but less well in the complex case. We conclude from our study that mixed-precision Runge--Kutta-type methods can yield a significant speedup without reducing the accuracy given that the matrix-vector products are carefully implemented to avoid unnecessary truncation errors and the solver tolerances are chosen above the machine epsilon of the chosen precision.

Future work entails extending to a more complicated setting, e.g., non-linear equations where the main challenge lies in deploying efficient block solvers for the implicit steps. This includes multigrid methods, incomplete factorizations, sparse pseudo-inverse method etc. The current setting of linear problems combined with FDM preconditioners is suitable for a well-controlled study demonstrating the importance of well optimized kernels to achieve good speedups, as we demonstrate in the performance analysis in Sections \ref{Sec:Tensor_comp} and \ref{Sec:speedupsolver}. Furthermore, in order to treat problems of practical interest, a distributed setting allowing for multiple GPUs or CPU nodes is also necessary which introduces communication overheads which may impact the potential speedup. Finally, investigation of integrating even lower precisions is of interest---here significant care has to be taken to avoid the stability issues outlined in \Cref{sec:stability} which will limit the amount of viable methods. Due to the necessity to prevent over-/underflow, lower precision approaches likely require additional checks introducing overhead and limiting the potential speedup.

\section*{Acknowledgement} 
 This work was funded by the German Ministry of Education and Research through the project ``PDExa: Optimized software methods for solving partial differential equations on exascale supercomputers'', grant agreement no.~16ME0639 and 16ME0640, and by the European Union---NextGenerationEU. Views and opinions expressed are however those of the author(s) only and do not necessarily reflect those of the European Union or the European Commission. Neither the European Union nor the European Commisson can be held responsible for them. The authors gratefully acknowledge the funding of this project by computing time provided by the Paderborn Center for Parallel Computing (PC2) through the project ``PDExa: Mixed precision'', project ID PDExaMP.

\bibliographystyle{abbrv}
\bibliography{bibfile}

\pagebreak
\begin{appendix}\section{Proof of A-stability of the midpoint rule with arbitrary number of correction steps}
\textit{A-stability is shown based on the stability function, e.g., the proof uses exact arithmetic. }
The implicit midpoint rule with $n-1$ corrector steps has the following coefficients
$$
A=\left(\begin{array}{ccccc}
\frac{1}{2}& 0 & 0& \hdots & 0\\
\frac{1}{2}& 0 & 0 & \hdots &0 \\
0&  \frac{1}{2} & 0 & \ddots & \vdots\\
\vdots & \ddots & \ddots & \ddots& 0\\
0& \hdots & 0 & \frac{1}{2} & 0
\end{array}\right)
\text{ and  }  
b=\left(\begin{array}{c}
0 \\
\vdots \\
0\\
1
\end{array}\right),
$$
where  $A \in \mathbb{R}^{n \times n}$ and $b \in \mathbb{R}^n$.  The stability function is given by $R(z)=1+z b^\top (I_n -zA)^{-1}e$, where $e=\left(\begin{array}{cccc}
1 & 1 & 1 & 1
\end{array}\right)^\top$ and $I_n$ denotes the identity matrix of size $n$. We have
$$
I_n -zA=\left(\begin{array}{ccccc}
\frac{2-z}{2}& 0 & 0& \hdots & 0\\
-\frac{z}{2}& 1 & 0 & \hdots &0 \\
0&  -\frac{z}{2} & 1 & \ddots & \vdots\\
\vdots & \ddots & \ddots & \ddots& 0\\
0& \hdots & 0 & -\frac{z}{2} & 1
\end{array}\right)
$$
The inverse $(I_n-zA)^{-1}$ can be directly calculated by by Gaussian elimination. Multiplying the first row by $\frac{-2}{z-2}$ and adding $\frac{z}{2}$ of the first row to the second yields:

\begin{align*}
&\left(\begin{array}{ccccc|ccccc}   
\frac{2-z}{2}& 0 & 0& \hdots & 0        & 1 & 0 & 0 & \hdots & 0   \\
-\frac{z}{2}& 1 & 0 & \hdots &0         & 0 & 1 & 0     &  &   \\
0&  -\frac{z}{2} & 1 & \ddots & \vdots  & \vdots & \ddots & \ddots & \ddots & \vdots  \\
\vdots & \ddots & \ddots & \ddots& 0    &         &  & & 1 & 0\\
0& \hdots & 0 & -\frac{z}{2} & 1        &    0    & \hdots & & 0 & 1 \\
\end{array}\right)
\Leftrightarrow \\ 
&\left(\begin{array}{ccccc|ccccc} 
1 & 0 & 0& \hdots & 0                   & \frac{-2}{z-2} & 0 & 0 & \hdots & 0   \\
0 & 1 & 0 & \hdots &0                   & \frac{-z}{z-2} & 1 & 0     &  &   \\
0&  -\frac{z}{2} & 1 & \ddots & \vdots  & \vdots & \ddots & \ddots & \ddots & \vdots  \\
\vdots & \ddots & \ddots & \ddots& 0    &         &  & & 1 & 0\\
0& \hdots & 0 & -\frac{z}{2} & 1        &    0    & \hdots & & 0 & 1 \\
\end{array}\right)
\end{align*}
The same two operations can be repeated on each row. The final form of the inverse can be found using the following identity; for any $k \in \mathbb{N}_{\geq0}$
$$
\frac{-2}{z-2} \cdot \left( \frac{z}{2} \right)^k = \frac{-z^k}{2^{k-1}\cdot z-2^k},
$$
which can be proven by means of induction: 
\begin{enumerate}
    \item $k=0$: $\frac{-2}{z-2} \cdot 1 =\frac{-z^0}{2^{-1}\cdot z-2^{-1}}$
    \item $k>0$: $\frac{-2}{z-2} \cdot \left( \frac{z}{2} \right)^{k+1} \overset{Ind.}{=} \frac{-z^k}{2^{k-1}\cdot z-2^k} \cdot \frac{z}{2}=\frac{-z^{k+1}}{2(2^{k-1}z-2^k)}=\frac{-z^{k+1}}{2^kz-2^{k+1}}$
\end{enumerate}
We thus have
$$
(I_n -z A)^{-1} =\left(\begin{array}{ccccc}
\frac{-2}{z-2} & 0 & 0 & \hdots & 0   \\
 \frac{-z}{z-2} & 1 & 0     &  &   \\
 \vdots & \vdots & \ddots & \ddots & \vdots  \\
\frac{-z^{n-2}}{2^{n-3}z-2^{n-2}}& \frac{z^{n-2}}{2^{n-2}} & \hdots & 1 & 0\\
 \frac{-z^{n-1}}{2^{n-2}z-2^{n-1}}& \frac{z^{n-1}}{2^{n-1}} & \hdots & \frac{z}{2} & 1\\
\end{array}\right).
$$
Next we show
\begin{equation}\label{veceq_Stab}
    (I_n-zA)^{-1}e = \left(\frac{-2}{z-2} \cdot \left( \frac{z}{2} \right)^k + \sum\limits_{i=0}^{k-1} \frac{z^{i}}{2^{i}}\right)_{k=0,\ldots,n-1}
=\frac{-2}{z-2}e
\end{equation}
by proving the following statement, again by the use of induction,
$$
\frac{-2}{z-2} \cdot \left( \frac{z}{2} \right)^k + \sum\limits_{i=0}^{k-1} \frac{z^{i}}{2^{i}} = \frac{-2}{z-2} \quad \forall k \in \mathbb{N}_{\geq0}.
$$
\begin{enumerate}
    \item $k=0$: $\frac{-2}{z-2} \cdot 1+0==\frac{-2}{z-2}$ 
    \item $k>0$: $\frac{-2}{z-2} \cdot \left( \frac{z}{2} \right)^{k+1}+\sum\limits_{i=0}^{k} \frac{z^{i}}{2^{i}}=\frac{-2}{z-2} \cdot \left( \frac{z}{2} \right)^{k+1}+\sum\limits_{i=1}^{k} \frac{z^{i}}{2^{i}} +1= \frac{z}{2} \Bigl( \frac{-2}{z-2} \cdot \left( \frac{z}{2} \right)^k + \sum\limits_{i=0}^{k-1} \frac{z^{i}}{2^{i}} \Bigr)+1 \overset{Ind.}{=} \frac{z}{2} \cdot \frac{-2}{z-2}+1=\frac{-2z}{2z-4}+\frac{2z-4}{2z-4}=\frac{-2}{z-2}$
\end{enumerate}
Eq. \eqref{veceq_Stab} directly follows. Further
$$
R(z)=1+z b^\top (I_n -zA)^{-1}e = 1+z \cdot \frac{-2}{z-2}=\frac{z-2}{z-2}+\frac{-2z}{z-2}=\frac{z+2}{2-z}
$$
A-stability requires $|R(z)| =\Big| \frac{z+2}{2-z} \Big|\leq 1 \, \ \forall z \text{ with } Re(z)\leq 0$. Let $z=x+iy$ we then have
\begin{eqnarray*}
    \frac{\sqrt{(x+2)^{2}+y^{2}}}{\sqrt{(2-x)^{2}+y^{2}}} \leq 1 &\Leftrightarrow& (x+2)^{2}+y^{2} \leq (2-x)^{2}+y^{2} \\
    & \Leftrightarrow & x^{2}+4x+4 \leq 4-4x+x^{2} \\
    & \Leftrightarrow & x \leq -x
\end{eqnarray*}
which holds true when $x \leq 0$, i.e., when $Re(z) \leq 0$. Thus A-stability is proven for any positive number of correction steps. 
$\square$  
\end{appendix}

\end{document}